\documentclass[10pt]{amsart}

\usepackage{amssymb,amsmath,amsthm}
\usepackage[all]{xy}  %for \xymatrix
\usepackage{url}
\usepackage{color}
\usepackage{mdframed} %shaded and framed boxes
\usepackage{fancybox} %Sbox
\usepackage{multirow}
\usepackage[table]{xcolor}
\usepackage{textcomp,listings}
\usepackage{versions}
\usepackage{rotating} %for includegraphics
\usepackage{placeins} %for \FloatBarrier
\usepackage{lineno} %\linenumbers
\usepackage{graphicx} %scalebox

\usepackage{algorithm}
\usepackage{algorithmicx, algpseudocode}

\definecolor{green}{rgb}{0,0.5,0}
\definecolor{dkgreen}{rgb}{0,0.6,0}
\definecolor{gray}{rgb}{0.5,0.5,0.5}
\definecolor{mauve}{rgb}{0.58,0,0.82}
\definecolor{questcolor}{rgb}{0.9, 0.9, 0.9}
\lstnewenvironment{python}[1][]{
\lstset{
  frame=single,                   % adds a frame around the code
  language=python,                          % the language of the code
  basicstyle=\ttfamily\small, % the size of the fonts that are used for the code
  numbers=left,                             % where to put the line-numbers
  numberstyle=\scriptsize\color{black},  % the style that is used for the line-numbers
  stepnumber=1,                   % the step between two line-numbers. If it's 1, each line
  numbersep=8pt,                  % how far the line-numbers are from the code
  backgroundcolor=\color{white},      % choose the background color.
  stringstyle=\color{red},
  showspaces=false,               % show spaces adding particular underscores
  showstringspaces=false,         % underline spaces within strings
  showtabs=false,                 % show tabs within strings adding particular underscores
  tabsize=2,                      % sets default tabsize to 2 spaces
  captionpos=b,                   % sets the caption-position to bottom
  breaklines=true,                % sets automatic line breaking
  breakatwhitespace=false,        % sets if automatic breaks should only happen at whitespace
  rulecolor=\color{black},        % if Ä set, the frame-color may be changed on line-breaks within not-black text (e.g. commens (green here))
  framexleftmargin=1mm, framextopmargin=1mm, frame=shadowbox, %rulesepcolor=\color{blue},#1
  alsoletter={1234567890},
  otherkeywords={\ , \}, \{},
  keywordstyle=\color{blue},       % keyword style
  stringstyle=\color{mauve},         % string literal style
  emph={access,and,break,class,continue,def,del,elif ,else,%
  except,exec,finally,for,from,global,if,import,in,i s,%
  lambda,not,or,pass,print,raise,return,try,while},
  emphstyle=\color{black}\bfseries,
  emph={[2]True, False, None, self},
  emphstyle=[2]\color{blue},
  emph={[3]from, import, as},
  emphstyle=[3]\color{blue},
  upquote=true,
  morecomment=[s]{"""}{"""},
  commentstyle=\color{dkgreen},       % comment style
  emph={[4]1, 2, 3, 4, 5, 6, 7, 8, 9, 0},
  emphstyle=[4]\color{blue},
  literate={:}{{\textcolor{blue}:}}{1}%
  {=}{{\textcolor{blue}=}}{1}%
  {-}{{\textcolor{blue}-}}{1}%
  {+}{{\textcolor{blue}+}}{1}%
  {}{{\textcolor{blue}}}{1}%
  {!}{{\textcolor{blue}!}}{1}%
  {(}{{\textcolor{blue}(}}{1}%
  {)}{{\textcolor{blue})}}{1}%
  {[}{{\textcolor{blue}[}}{1}%
  {]}{{\textcolor{blue}]}}{1}%
  {<}{{\textcolor{blue}<}}{1}%
  {>}{{\textcolor{blue}>}}{1},%
}}{}

\newcommand{\hhat}{{\hat h}}

\providecommand{\abs}[1]{\left\lvert#1\right\rvert}
\providecommand{\norm}[1]{\left\lvert\left\lvert#1\right\rvert\right\lvert}

\usepackage{listofitems}
\newcommand{\sm}[1]{%
  \readlist\smvar{#1}%
  \left(\begin{smallmatrix}\smvar[1]&\smvar[2]\\ \smvar[3]&\smvar[4]\end{smallmatrix}\right)}

\def\C{\mathbb{C}}

\def\P{\mathbb{P}}
\def\Q{\mathbb{Q}}

\def\calM{\mathcal{M}}

\def\calO{\mathcal{O}}

\def\Kbar{\overline{K}}

\DeclareMathOperator{\Fix}{Fix}

\DeclareMathOperator{\PGL}{PGL}
\DeclareMathOperator{\Pre}{Pre}

\theoremstyle{plain}% default
\newtheorem{thm}{Theorem}

\newtheorem*{openquestion}{Open Questions}
\theoremstyle{definition}

\newtheorem{conj}{Conjecture}
\newtheorem{exmp}[thm]{Example}

\theoremstyle{remark}

\title[Extreme Examples in Arithmetic Dynamics]{A Genetic Algorithm for Generating Extreme Examples in Arithmetic Dynamics}
\author{Benjamin Hutz}
\address{Department of Mathematics and Statistics,
         Saint Louis University,
         St.~Louis, MO 63103}
\email{benjamin.hutz@slu.edu}

\includeversion{all_data}
\begin{document}

\keywords{periodic points, canonical height, genetic algorithm, extreme examples, dynamical system, preperiodic points}

\subjclass[2020]{
37-04,
    % Software, source code, etc. for problems pertaining to dynamical systems and ergodic theory 
 37P15,
    % Dynamical systems over global ground fields
 37-11,
    % Research data for problems pertaining to dynamical systems and ergodic theory
}

%\date{\today}

    \begin{abstract}
        We describe a genetic algorithm to find extreme examples in the arithmetic of dynamical systems. The algorithm is applied to four problems: small (non-zero) canonical heights, many rational preperiodic points, long rational cycles, and long rational tails. Data is provided for extreme examples generated for polynomials up to degree 13 and rational functions up to degree 5. This work significantly expands the known examples of extreme behavior for several of the conjectured behaviors in arithmetic dynamics and provides a foundation from which to begin a more advanced application of machine learning techniques in the creation of extreme examples for arithmetic dynamics. Further, the algorithm can be easily applied to other problems by replacing the scoring function.
    \end{abstract}

\maketitle

With the rapid improvements in machine learning and some recent success with applying neural networks to extreme examples in graph theory \cite{GAKS, Wagner}, we examine the feasibility of generating extreme examples in arithmetic dynamics\footnote{Thanks to Trevor Hyde for posing the problem to me and making me aware of the recent results in graph theory.}. The field of arithmetic dynamics is concerned with iterating polynomial or rational maps. Let $K$ be a field and $f \in K(x)$, be a polynomial or rational function. We define the $n$th iterate of $f$ as $f^n = f \circ f^{n-1}$ for $n> 0$ and $f^0(x) = x$. A point $P$ is \emph{preperiodic} if the forward orbit $\calO = \{f^n(P) : n \geq 0\}$ is finite and \emph{wandering} otherwise. For brevity we will call such a function $f$ a \emph{dynamical system}. Silverman's graduate textbook \cite{Silverman10} gives a good introduction to the field. There are a number of broad open questions in the arithmetic of dynamical systems, such as the Uniform Boundedness Conjecture of Morton and Silverman \cite{Silverman7} or a dynamical formulation of Lang's conjecture on canonical heights \cite[Conjecture 4.98]{Silverman10}. These problems will be discussed in detail in Section \ref{sect_problems}. Both of these conjectures, among others, propose the existence of certain constants limiting the extreme behavior of some arithmetic property of dynamical systems. For example, the Morton-Silverman conjecture limits the total number of rational preperiodic points of a dynamical system and the dynamical version of Lang's conjecture limits the minimum canonical height of a wandering point. The constants in these conjectures depend, of course, on various intrinsic properties of the dynamical system, such as the degree. For small degrees (2 or 3) various authors have computationally examined these conjectures, mainly through extensive searches, e.g., \cite{Benedetto8, Benedetto4, DFK, Hutz5, Poonen}. However, as the degree increases, this approach quickly becomes impractical. Getting a better understanding of the constants in these conjectures through extreme examples would be beneficial to the field. For example, Poonen's classification of all $\Q$-rational preperiodic structures for quadratic polynomials \cite{Poonen} (similar to Mazur's theorem for $\Q$-rational torsion on elliptic curves \cite{Mazur}) depends on the conjecture that there are no $\Q$-rational periodic points with cycle length longer than 3 for a quadratic polynomial $f \in \Q[x]$. Our results show that even a relatively simple genetic algorithm is successful in creating extreme examples up to moderate degree. We examine extreme behavior for four problems in the arithmetic of dynamical systems:
\begin{enumerate}
    \item Small canonical height for wandering points
    \item Many rational preperiodic points
    \item Long rational periodic cycles
    \item Long rational preperiodic tails
\end{enumerate}
for
\begin{itemize}
    \item polynomials up to degree $13$
    \item rational functions up to degree $5$.
\end{itemize}
Note that the largest successful degree depended somewhat on the problem, so extreme examples for the highest degree were not achieved for all four problems. Section \ref{sect_data} summarizes the most extreme examples found, including a comparison to the best previously known (as far as the author is aware) for each problem in each degree. We summarize the trend data on the open conjectures in the following general question split into subquestion for each problem and the rational versus polynomial cases. Question (2)(a) comes from Doyle-Hyde \cite{Doyle-Hyde}, and by extension so does (2)(b), the remaining are derived from the data obtained from this algorithm. All these questions can be used to gauge efficacy of future algorithms on these problems.
\begin{openquestion}
    Prove the following general bounds or find a more extreme example.
    \begin{enumerate}
        \item \textbf{Small canonical heights}
        \begin{enumerate}
            \item For a degree $d\geq 2$ polynomial map defined over $\Q$, the smallest height ratio grows asymptotically as $10^{-\frac{3d}{2}}$.
            \item For a degree $d\geq 2$ rational map defined over $\Q$, smallest height ratio grows asymptotically as $10^{-2d}$ for rational maps.
        \end{enumerate}
        \item \textbf{Many rational preperiodic points}
        \begin{enumerate}
            \item A degree $d\geq 2$ polynomial map defined over $\Q$ exists with at least $d + \max(6, \lfloor \log_2(d)\rfloor)$ $\Q$-rational periodic points.
            \item A degree $d\geq 2$ rational map defined over $\Q$ exists with at least $d + \max(6, \lfloor \log_2(d)\rfloor)$ $\Q$-rational periodic points.
        \end{enumerate}
        \item \textbf{Long rational periodic cycles}
        \begin{enumerate}
            \item A degree $d\geq 2$ polynomial map defined over $\Q$ has no $\Q$-rational preperiodic point with minimal period larger than $d+3$. 
            \item A degree $d\geq 2$ rational map defined over $\Q$ has no $\Q$-rational preperiodic point with minimal period larger than $2d+3$.
        \end{enumerate}
        \item \textbf{Long rational preperiodic tails}
        \begin{enumerate}
            \item A degree $d\geq 2$ polynomial map there were no examples of a $\Q$-rational preperiodic point with preperiodic tail larger than $d+2$. 
            \item A degree $d\geq 2$ rational map there were no examples of a $\Q$-rational preperiodic point with preperiodic tail larger than $2d+2$.
        \end{enumerate}
    \end{enumerate}
\end{openquestion}

The basic idea of the algorithm is the following general steps of a genetic algorithm:
\begin{enumerate}
    \item Generate a random set of dynamical systems
    \item Score those dynamical systems from best to worst
    \item Keep only the best scoring systems
    \item Create the next generation by combining those best systems to create new (and hopefully better) examples.
\end{enumerate}
A principle tenet for this project was to keep the algorithm as simple as possible and to only change the ``scoring'' function as we examine the four different problems. This allows for a flexible algorithm that could be applied by other researchers to generate their own extreme examples by either using our scoring functions, or creating their own to plug into the general algorithm.

The article is organized as follows. Section \ref{sect_problems} describes the mathematical background on each open problem in arithmetic dynamics that we consider. A discussion of the current state of knowledge is included. Section \ref{sect_algorithm} describes the details of the algorithm. Section \ref{sect_data} provides a summary of the most extreme examples found with this algorithm.
A fuller list of examples is provided as an appendix to the arxiv version\footnote{arXiv:2601.11482}.
The code for the algorithm as implemented for SageMath \cite{sage} can be found on github\footnote{\url{https://github.com/bhutz/GA_for_ADS_examples}}.

\section{Open Problems to Arithmetic Dynamics} \label{sect_problems}

In this section, we give an overview of the problems to be examined including their current status. A common theme is that these problems conjecture the existence of a constant bounding some arithmetic property of dynamical systems. With extensive data only for degree 2 maps and degree 3 polynomials, the precise behavior of these constants with respect to their dependencies (such as degree) remains unknown.

\subsection{Small Height Ratio}

        Let $K$ be a number field and recall that the height of a point $P \in \P^N(K)$ is defined as
        \begin{equation*}
            H_K(P) = \prod_{v \in M_k}\norm{P}_v^{n_v} = \prod_{v \in M_k}\max(\abs{P_0}_v, \cdots, \abs{P_N}_v)^{n_v},
        \end{equation*}
        where $M_K$ is the set of all absolute values and $n_v$ is the local degree. We set $h_K(P) = \log(H_K(P))$. This height can be normalized by the degree of the field as $H(P) = H_K(P)^{1/[K:\Q]}$ to produce an absolute height. 
        The height $H(P)$ is reasonably well behaved under iteration by morphisms $f:\P^N(K) \to \P^N(K)$ of degree $d \geq 2$ in the sense that there exists a constant $C_1$ that depends on $f$, but not $P$, so that
        \begin{equation*}
            \abs{h(f(P)) - dh(P)} < C_1, \quad \forall P \in \P^N(K).
        \end{equation*}
        In many applications we would prefer an exact functional equation, so we define the canonical height as
        \begin{equation*}
            \hhat_f(P) = \lim_{n \to \infty} \frac{h(f^n(P))}{\deg(f)^n}.
        \end{equation*}
        The canonical height satisfies the exact functional equation
        \begin{equation*}
            \hhat_f(f(P)) = d \hhat_f(P), \quad \forall P \in \P^N(K).
        \end{equation*}
        This functional equation immediately implies that preperiodic points have canonical height 0. In the case of global fields, the reverse implication is also true. Furthermore, 
        there exists a constant $C_2$ depending on $f$, but not $P$ so that
        \begin{equation} \label{eq2}
            \abs{\hhat_f(P) - h(P)} < C_2, \quad \forall P \in \P^N(K).
        \end{equation}
        See Silverman \cite[Chapter 3]{Silverman10} for an introduction to heights and canonical heights for dynamical systems.
        Since the set of points of $\P^N(K)$ of bounded height is finite, Equation \eqref{eq2} implies there must exist a wandering point $P$ with smallest positive canonical height. Consequently, an interesting question is how small can this smallest positive canonical height be. This is similar to Lang's conjecture for elliptic curves which conjectures, for any non-torsion point $P$ and an elliptic curve defined over a number field $K$, the existence of constants $C_3,C_4>0$ so that
        \begin{equation*}
            \hhat(P) > C_3 \log(\mathcal{N}_K(\mathcal{D}_E)) + C_4,
        \end{equation*}
        where $\mathcal{D}_E$ is the minimal discriminant of the elliptic curve $E$. In the dynamical system case, some additional care needs to be taken since conjugating the function can lead to arbitrarily small canonical heights. Recall that we say two dynamical systems $f,g: \P^N(K) \to \P^N(K)$ are \emph{conjugate} if there is an element $\alpha \in \PGL_{N+1}(\Kbar)$ so that 
        \begin{equation*}
            f^{\alpha} = \alpha^{-1} \circ f \circ \alpha = g.
        \end{equation*}
        We define the moduli space of dynamical systems of degree $d$, $\calM_d$, as the set of conjugacy classes under this conjugation action \cite{Levy, petsche}.
        The relevant conjecture for dynamical systems is as follows
     \begin{conj}[\textup{\cite[Conjecture 4.98]{Silverman10}}] \label{conj-lang}
        For any $d \geq 2$, there is a constant $M(d)>0$ depending only on $d$ so that for any $f:\P^1 \to \P^1$ of degree $d$ and wandering point $Q$ we have
        \begin{equation} \label{eq1}
            \hhat_f(Q) \geq M(d) h_{\mathcal{M}}(f)
        \end{equation}
        where $h_{\mathcal{M}}(\cdot)$ is a height on the moduli space of dynamical systems of degree $d$.
    \end{conj}
    Roughly speaking, a height on the moduli space is a height function that is invariant on conjugacy classes. For a more detailed description, see Silverman \cite{Silverman20}. There are a number of different choices for heights on moduli space. For example, the space of all quadratic polynomials is parameterized as $f_c(z) = z^2+c$, and this moduli space is isomorphic to $\C$. In this case, $h_{\mathcal{M}}(f_c) = h(c)$ defines a height on moduli space. Alternatively, we can consider the multipliers. We can define the set of multipliers for the fixed points of $f$ as $\{ f'(P) : P \in \Fix(f)\}$. If we then take the elementary symmetric functions applied to these multipliers, we end up with invariants of the moduli space. For degree $2$ maps, the first two fixed point multiplier invariants $(\sigma_1,\sigma_2)$ were shown by Milnor \cite{Milnor} to give an isomorphism $\C^2 \cong \mathcal{M}_2$. Consequently, a height on moduli space is given by $h_{\mathcal{M}}(f) =  \max(h(\sigma_1),h(\sigma_2))$. More generally, McMullen proved that these invariants, taken from periodic points of large enough period, determine a conjugacy class up to finitely many choices (outside of the Latt\`es maps) \cite{McMullen2}. Consequently, we can use these multiplier invariants via explicit forms of McMullen's theorem to formulate a general height on moduli space. For example, for polynomials, the fixed point multiplier invariants determine the conjugacy class up to finitely many choices (Fujimura and Nishizawa \cite{Fujimura3}). For rational functions of degree greater than $2$, the fixed point multiplier invariants, no longer determine finitely many conjugacy classes, however, the first and second periodic point multiplier invariants determine finitely many conjugacy classes: Ji-Xie conjectured in 2023 that the first and second multiplier invariants should be generically one-to-one \cite{Ji}. To construct the height function more formally, as in Ingram \cite{Ingram6} or Benedetto-Ingram-Jones-Levy \cite{Benedetto5}, consider the morphism
    \begin{align*}
        \Lambda_2: \calM_d &\to \P^{d+1} \times \P^{d^2+1}\\
        \Lambda_2(f) &= (\boldsymbol \sigma^{(1)}, \boldsymbol \sigma^{(2)}),
    \end{align*}
    where $\boldsymbol \sigma^{(n)}$ are (the homogenization of) the elementary symmetric functions applied to the multipliers of the points of period $n$ for $f$. As just discussed, this map is finite (outside of the Latt\`es locus) via a refinement of McMullen's theorem and is expected to be injective (outside of the Latt\`es locus). Consequently, we can pullback the regular height on $\P^{d+1} \times \P^{d^2+1}$ to get a height on moduli space (outside of the Latt\`es locus). For simplicity, and to compare the polynomial and rational cases, we continue to choose $h_{\mathcal{M}}$ as the height of just the fixed point multiplier invariants even though this is not a true height on moduli space for general rational functions.

    To examine the relation expected in the dynamical version of Lang's conjecture (Conjecture \ref{conj-lang}), it makes sense to look at the smallest possible value of the ratio of heights
    \begin{equation*}
        \frac{\hhat_f(P)}{h_{\mathcal{M}}(f)}.
    \end{equation*}
    This ratio has been examined in the cases of quadratic polynomials \cite{Benedetto8}, cubic polynomials \cite{Benedetto4}, and quadratic rational functions \cite{Benedetto8}. The choice of $h_{\mathcal{M}}$ varied among the different investigations so can affect the specific numerical values. To compare our results numerically to previous results, we recompute all minimum values with our $h_{\mathcal{M}}$.

    \subsection{Many Rational Preperiodic Points (Morton-Silverman Conjecture)}

    In 1994 Morton and Silverman formulated a very general conjecture similar to the boundedness of the number of rational torsion points on elliptic curves. Specifically, they conjectured that the number of rational preperiodic points should be uniformly bounded with the constant depending only on the degree and dimension.
    \begin{conj}[\cite{Silverman7}]
        For any integers $d \geq 2$ and $D \geq 1$, there exists a constant $C(d,D,N)$ so that for any endomorphism $f:\P^N \to \P^N$ of degree $d$ we have 
        \begin{equation*}
            \#\Pre(f,K) \leq C(d,D,N)
        \end{equation*}
        for any number field $K$ with $[K:\Q] \leq D$.
    \end{conj}

    To get an idea of the scope of this conjecture, note that the existence of a constant $C(4,1,1)$ implies Mazur's theorem on rational torsion of elliptic curves (without the specific classification of subgroups) and $C(4,D,1)$ implies Merel's theorem on the uniform boundedness of rational torsion for elliptic curves.
    
    We re-state the Morton-Silverman conjecture on the uniform boundedness of rational preperiodic points in the simple case of the rational numbers and dimension 1.
    \begin{conj}
        For any integer $d \geq 2$, there exists a constant $C(d)$ depending only on $d$ so that
        for any endomorphism $f:\P^1 \to \P^1$ of degree $d$ we have 
        \begin{equation*}
            \#\Pre(f,\Q) \leq C(d).
        \end{equation*}
    \end{conj}
    Even in the special case of quadratic polynomials, this conjecture remains open. Poonen demonstrated that if there are no quadratic polynomials with rational periodic points with minimal period greater than $3$, then $C_{poly}(2) = 9$ \cite{Poonen}. It was shown by Morton \cite{Morton2} that there are no quadratic polynomials with a $\Q$-rational point of minimal period $4$ and by Flynn-Poonen-Schaefer that there are none with a period 5 point \cite{FPS}, and conditionally by Stoll that there are none with a period $6$ point  \cite{Stoll3}.

    Quadratic rational functions have been explored by a few authors. Manes in her ph.d. thesis conjectures that $C(2) \leq 12$ with methods similar to Poonen \cite{Manes2}. Benedetto et al \cite{Benedetto8} has perhaps the most complete set of data giving maps with up to 14 $\Q$-rational preperiodic points and cycles of length up to $7$. Blanc-Canci-Elkies give an infinite family with 14 $\Q$-rational preperiodic points \cite{BCE}.

    Benedetto et al \cite{Benedetto8} parameterized the maps by examining orbits which begin with [$\infty \to 1 \to 0 \to x_3 \to x_4 \to x_5$] with $H(x_3,x_4,x_5) \leq 100$. Each such orbit uniquely determines a quadratic map. There are approximately 1.7 trillion such triples making this an extensive search.

    Benedetto et al \cite{Benedetto4} examined cubic polynomials and found at most 11 $\Q$-rational preperiodic points. They looked at $az^3+bz+1$ or $az^2+bz$ for $H(a,b) \leq 300$. There are approximately 24 billion such maps making this also an extensive search.

    Besides this work to compute a value for specific $C(d,D,N)$ in various special cases, there has also been some effort to determine a minimum growth rate for $C(d,D,N)$ in terms of the dependencies $d$, $D$, and $N$. We enumerate those efforts here.
    \begin{enumerate}
        \item $C(d,D,N)$ grows at least linearly in $D$ \cite{Hutz5}.
        \item $C(d,D,N)$ grows at least linearly in $d$ by interpolation and at least as $d + \log_2(d)$ for $d$ large enough by \cite{Doyle-Hyde}.
        \item $C(d,D,N)$ grows faster than $C(k)N^k$ for some constant $C$ depending on $k$ for any integer $k$ \cite{Hutz4}.
    \end{enumerate}
    Most related to the Morton-Silverman conjecture and this work is the dependency on $d$. For $d\geq 2$, Doyle and Hyde \cite{Doyle-Hyde} define
    \[
        B_d := \sup_f |\Pre(f,\Q)|
    \]
    where the supremum is taken over all polynomials $f(x) \in \Q[x]$ with $2 \leq \deg(f) \leq d$.
    They prove that for all $d \geq 2$ we have
    \[
        B_d \geq d + \max(6, \lfloor \log_2(d)\rfloor).
    \]
    By construction, $B_d$ is the optimal value of $C(d,1,1)$ restricted to polynomials.
    A simple interpolation argument implies that $B_d \geq d + 1$; Doyle-Hyde show that $B_d - d$ is unbounded and grows at least as fast as $\log(d)$.
    Their methods were not provably optimal. Hence, we would like to
    \begin{enumerate}
        \item Find explicit examples with an extreme number of rational preperiodic points, and
        \item Gather data on the true lower bound for $B_d$.
    \end{enumerate}
    In particular, the lower bound they give is simply $B_d \geq d + 6$ for $d < 128$, which comes from a single family of examples. We were able to equal this lower bound for $2 \leq d \leq 5$ and did not find any examples that improve upon this bound.

    \subsection{Long Periodic Cycle}

        A key step in the known partial classifications of rational preperiodic structures related to the Morton-Silverman uniform boundedness conjecture is bounding the size of a rational cycle \cite{Manes2, Poonen}. Furthermore, there is known to be a (non-explicit) uniform bound on the number of rational preimages \cite{FHIJMTZ}, so a uniform bound on rational cycle length would be a strong step towards resolving the Morton-Silverman conjecture. Consequently, it is quite interesting to look for examples of maps with long cycles of rational points. Interpolation guarentees a cycle length of at least $d+1$ for polynomials and $2d+2$ for rational functions. So the question is concerning how much larger is achievable than this trivial bound. The previous computational investigations on quadratic and cubic polynomials and quadratic rational maps includes data on long rational cycles and is included in the summary data tables in Section \ref{sect_data}.
    
    \subsection{Long Preperiodic Tail}

        The result that the length of a rational periodic tail is uniformly bounded is not explicit \cite{FHIJMTZ}. Finding an upper bound on the tail length when combined with an explicit bound on the cycle length is one way to give an explicit upper bound on the constant $C(d,D,N)$ for Morton-Silverman. Consequently, we look for examples with very long rational tails. The computational investigation on quadratic and cubic polynomials and quadratic rational maps included data on long rational tails and is included in the summary data tables in Section \ref{sect_data}.

\section{The Algorithm} \label{sect_algorithm}
This section presents details of our algorithm. The method is a general genetic algorithm with fitness functions modified for each subproblem. Full code is available on github\footnote{\url{https://github.com/bhutz/GA_for_ADS_examples}} and can be used with SageMath.

Generally speaking, a genetic algorithm follows the following outline
\begin{enumerate}
    \item Start with a population of potential examples
    \item Rank the population in terms of a fitness function
    \item Keep a small fraction of the best performing examples: the survivors
    \item Repopulate by recombining the survivors with a small percentage of random changes (mutations)
    \item repeat
\end{enumerate}
We also found it helpful to periodically keep only the very best and regenerate the rest of the population randomly

\subsection{Population: Representation of Functions}
We chose to represent our functions via their orbits (similar to Benedetto et al \cite{Benedetto8}). Each such orbit, of appropriate length, corresponds to a unique polynomial or rational function through interpolation. Specifically, given a rational function with symbolic coefficients, every image $f(z_i) = z_{i+1}$ defines a linear equation in the coefficients. With enough independent linear equations, a function is uniquely determined. In the case of polynomials, this is the familiar method of Lagrange interpolation. A degree $d$ polynomial typically requires $d+1$ images and a degree $d$ rational function requires $2d+2$ images. We assume that all orbits start with $0$ and the first $d+1$ or $2d+2$ elements are integers. These assumptions are not restrictive, since any orbit can be conjugated to start at $0$ and further conjugated to clear denominators in the first finitely many elements. So our population consists of tuples of integers.

A disadvantage of this choice is that the same conjugacy class can have multiple representatives in our population. Since we are simply looking for the best possible examples, this did not pose any difficulties and our data set excludes these repetitions of conjugate functions.

\subsection{Recombination}
    An important step in any genetic algorithm is recombining the survivors to fill out the new generation. We explored two different methods of recombination
    \begin{enumerate}
        \item crossover
        \item permutation
    \end{enumerate}

    With crossover, we take two orbits and create a new orbit by taking the first half of one combined with the second half of the other.
    \begin{exmp}
        Given the orbits $O_1=[0,3,5,4,6,2,1]$ and $O_2 = [0,7,2,-1,9,8,3]$ which define degree $5$ polynomials we can recombine them in the two different methods. For crossover we take
        \begin{align*}
            [0,3,5,4,6,2,1] &\longrightarrow [3,5,4] \cup [6,2,1]\\
            [0,7,4,-1,9,8,3] &\longrightarrow [7,4,-1] \cup  [9,8,3].
        \end{align*}
        We then get the two new orbits
        \begin{align*}
            [0,3,5,4,9,8,3]\\
            [0,7,4,-1,6,2,1]
        \end{align*}
    \end{exmp}
    Alternatively, we can simply take a random permutation of the concatenation of the two orbits.
    \begin{exmp}
        With the same orbits $\calO_1$ and $\calO_2$ via permutation: we take the $12$ element list
        \begin{equation*}
            [3,5,4,6,2,1,7,4,-1,9,8,3]
        \end{equation*}
        and randomly permute to get
        \begin{equation*}
            [6,7,9,3,4,2,  -1,5,1,3,8,4]
        \end{equation*}
        which we split into the two new orbits
        \begin{align*}
            &[0,6,7,9,3,4,2]\\
            &[0,-1,5,1,3,8,4].
        \end{align*}
    \end{exmp}

    Furthermore, a fixed percentage of the time, elements of the orbit are replaced with a random integer to further assist with escaping from local extrema. We did not see any difference in results between the two recombination methods. However, the permutation method was slightly more costly from a performance perspective for larger degrees.

\subsection{Fitness Functions}
A key part of a genetic algorithm is the fitness function. This is the function that ranks the members of the current population. Choosing a fitness function for each problem that allowed successive generations to converge to an extreme example took some care. In the three preperiodic problems, the expectation is that a randomly generated map will have a wandering orbit for 0. So we needed a fitness function that would in some sense prefer functions where $0$ was ``almost'' preperiodic. The fitness function for small canonical heights was simpler in that you can simply compute the canonical height of $0$ for any randomly generated map and favor those with small (nonzero) values. We precisely describe our fitness functions in this section.

\subsubsection{Small height ratio}
We are looking for points with small canonical height ratio, and every orbit corresponding to a function has a computable canonical height for $0$. In the case that $0$ is preperiodic, this is undesirable for this problem, so we assign an arbitrary large score to those orbits. In the case that $0$ is wandering, we can simply compute the height ratio and assign the score to be that value. In this way, a lower score is better.

\begin{exmp}
 Consider a degree 4 polynomial. The orbit
    \begin{equation*}
        \calO = [3, 4, 5, 1, -1]
    \end{equation*} 
    corresponds to the map
    \begin{equation*}
        (x^4 - 28x^3y + 164x^2y^2 - 257xy^3 + 90y^4 : 30y^4)
    \end{equation*}
    The first multiplier invariants are
    \begin{equation*}
        [200, -78823/90, -105643/45, 2505080863/270000, 0]
    \end{equation*}
    so we take
    \begin{equation*}
        h_M(f) = \log(2505080863) = 21.6416.
    \end{equation*}
    The canonical height of $0$ is computed as
    \begin{equation*}
        \hhat_f(0) = 0.0002875
    \end{equation*}
    so the score of this function is the ratio
    \begin{equation*}
         \phi(\mathcal{O}) = \frac{\hhat_f(0)}{h_m(f)} \approx \frac{0.0002875}{21.6416} \approx 0.00001328.
    \end{equation*}
\end{exmp}

\subsubsection{Morton-Silverman Conjecture}
 The fitness function for the many rational preperiodic points problem covers two cases. If the orbit we are examining is preperiodic, then we can assign a score based on the total number of rational preperiodic points. If the orbit is wandering, then we need a different metric. In this case, we chose an arbitrary element in the orbit (past the interpolation limit) and compute its global height. The idea is that if the elements of the orbit remain small, then the orbit more likely to be preperiodic. So the fitness function is defined as follows
\begin{equation*}
    \phi(\mathcal{O}) = \begin{cases}
        h(f^n(0)) & $0$ \text{ is wandering and $n > d+1$ (or $(2d+2)$ for rational)}\\
        -\#\Pre(f,\Q) & \text{$0$ is a preperiodic point.} 
    \end{cases}
\end{equation*}

\begin{exmp}
 Consider a degree 4 polynomial. The orbit
    \begin{equation*}
        \mathcal{O} = [-3, -1, -5, -2, -4]
    \end{equation*} 
    corresponds to the map
    \begin{equation*}
        (-9x^4 - 34x^3y + 141x^2y^2 + 406xy^3 - 360y^4 : 120y^4)
    \end{equation*}
    $O$ is not preperiodic, so we return the height of the $6$th iterate
    \begin{equation*}
        \phi(\mathcal{O}) = h(f^6(0)) \approx 1.7918.
    \end{equation*}
    However the orbit
    \begin{equation*}
        \mathcal{O} = [-4, 1, -3, -1, -5]
    \end{equation*}
    corresponds to the map
    \begin{equation*}
        (-2x^4 - 11x^3y + 2x^2y^2 + 31xy^3 - 80y^4 : 20y^4)
    \end{equation*}
    The point $0$ is preperiodic for this map with period $(3,3)$. We compute, as via \cite{Hutz12}, the complete set of $\Q$-rational preperiodic points and find that there are a total of $7$, so we return a score of $\phi(\mathcal{O}) = -7$ for this map.
\end{exmp}

\subsubsection{Long Periodic Cycles}
The fitness function for the long periodic cycle problem is similar to the Morton-Silverman problem. If the orbit we are examining is preperiodic, then we can assign a score based on the length of the cycle. If the orbit is wandering, then we assign the fitness function value to be the global height of an arbitrary element in the orbit (past the interpolation limit). So the fitness function is defined as follows
\begin{equation*}
    \phi(\mathcal{O}) = \begin{cases}
        h(f^n(0)) & $0$ \text{ is wandering and $n > d+1$ (or $(2d+2)$ for rational)}\\
        -5n - m & \text{$0$ is a preperiodic point with minimal period $(m,n)$.} 
    \end{cases}
\end{equation*}

\begin{exmp}
    Consider a degree 7 polynomial. The orbit
    \begin{equation*}
        \mathcal{O} = [-5, 3, -1, -4, 2, 4, -3, 1]
    \end{equation*} 
    corresponds to the map
    \begin{align*}
        (47x^7 &+ 310x^6y - 934x^5y^2 - 8540x^4y^3 - 97x^3y^4 \\
        &+ 62950x^2y^5 + 35544xy^6 - 100800y^7 : 20160y^7)
    \end{align*}
    The orbit of $0$ is wandering and we compute the global height of the 9th element in the orbit as 
    \begin{equation*}
       h(f^9(0)) = 1.946.
    \end{equation*}
    We take $\phi(\mathcal{O}) =  1.946$.
    
    The orbit
    \begin{equation*}
        \mathcal{O} = [-5, 1, -1, -4, 2, 3, 4, -2]
    \end{equation*}
    corresponds to the map
    \begin{align*}
        (41x^7 &- 35x^6y - 1477x^5y^2 + 3535x^4y^3 + 16604x^3y^4\\
         &- 53900x^2y^5 - 45408xy^6 + 100800y^7 : -20160y^7)
    \end{align*}
    Now $0$ is preperiodic with period $(6,3)$ so the score for this map is now $\phi(\mathcal{O}) = -21$.
\end{exmp}

\subsubsection{Long Preperiodic Tails}
The fitness function for the long preperiodic tail problem is very similar to the periodic cycle problem. The difference is only in how the preperiodic case is scored with the emphasis on the tail rather than the periodic component. The fitness function is defined as follows
\begin{equation*}
    \phi(\mathcal{O}) = \begin{cases}
        h(f^n(0)) & $0$ \text{ is wandering and $n > d+1$ (or $(2d+2)$ for rational)}\\
        -n - 5m & \text{$0$ is a preperiodic point with minimal period $(m,n)$.} 
    \end{cases}
\end{equation*}

\begin{exmp}
        Consider a degree 7 polynomial. The orbit
    \begin{equation*}
        \mathcal{O} = [3, -3, 5, -1, 6, 1, -2, 4]
    \end{equation*} 
    corresponds to the map
    \begin{align*}
        (x^7 &+ 546x^6y - 5474x^5y^2 + 3360x^4y^3 + 84889x^3y^4\\
        &- 94626x^2y^5 - 442296xy^6 + 272160y^7 : 90720y^7)
    \end{align*}
    The orbit of $0$ is wandering and we compute the global height of the 9th element in the orbit as 
    \begin{equation*}
       h(f^9(0)) = 2.1972.
    \end{equation*}
    We take $\phi(\mathcal{O}) = 2.1972$.
    
    The orbit
    \begin{equation*}
        \mathcal{O} = [6, 1, 7, -1, 4, 8, -2, 5]
    \end{equation*}
    corresponds to the map
    \begin{align*}
        (-5x^7 &+ 105x^6y - 719x^5y^2 + 1515x^4y^3 + 1060x^3y^4\\
        &- 6660x^2y^5 + 14784xy^6 + 60480y^7 : 10080y^7)
    \end{align*}
    Now the $9$th element of the orbit is $5$ which is a fixed point so that $0$ is preperiodic with period $(8,1)$, so the score for this map is now $\phi(\mathcal{O}) = -41$.
\end{exmp}

\subsection{Discussion}

    This section contains a few observations about the performance of the algorithm. The computations were done on the Saint Louis University High Performance Cluster. Two types of jobs were run. Either several independent jobs with a population of 500 with 500 generations or a single jobs with a population of 2000 with 2000 generations. They were both successful at finding extreme examples with the larger populations requiring significantly more resources in terms of RAM and execution time. We had some preference for the small populations with fewer generations as convergence to an example was typically fast and it helped to conserve compute resources in the higher degrees. We summarize our main experience in the following points.
    \begin{itemize}
        \item A given population seemed to converge fairly quickly to a local best element. So it was helpful to periodically reset the population by keeping only the top few performers and randomly regenerating the rest of the population.
        \item Except for the Uniform Boundedness Problem, the algorithm appears to be finding the optimal or near optimal values based on comparison to the best currently known.
        \item The best known in the Uniform Boundedness problem appear to be maps with fairly dense preperiodic graphs. The algorithm mainly found either long orbits with few preimages, or long tails with few preimages. These do not appear to lead to optimal Uniform Boundedness values.
        \item When the number of parameters in the space reached around 12 (degree 11 polynomials or degree 5 rational functions) or higher, the algorithm had difficulty converging to an optimal value and mainly ended up getting stuck in local minimums.
    \end{itemize}

    As a comparison, we include the following two graphs which show a score comparison between the genetic algorithm and generating the same number of random maps. These were taken for degree 7 polynomials for the long periodic cycle problem and the small height ratio problem. In both cases, the genetic algorithm is significantly outperforming random sampling. Both of these graphs are typical of the comparisons we ran. This typicalness includes the very quick convergence of the periodic problem.

\begin{center}
\includegraphics[scale=0.45]{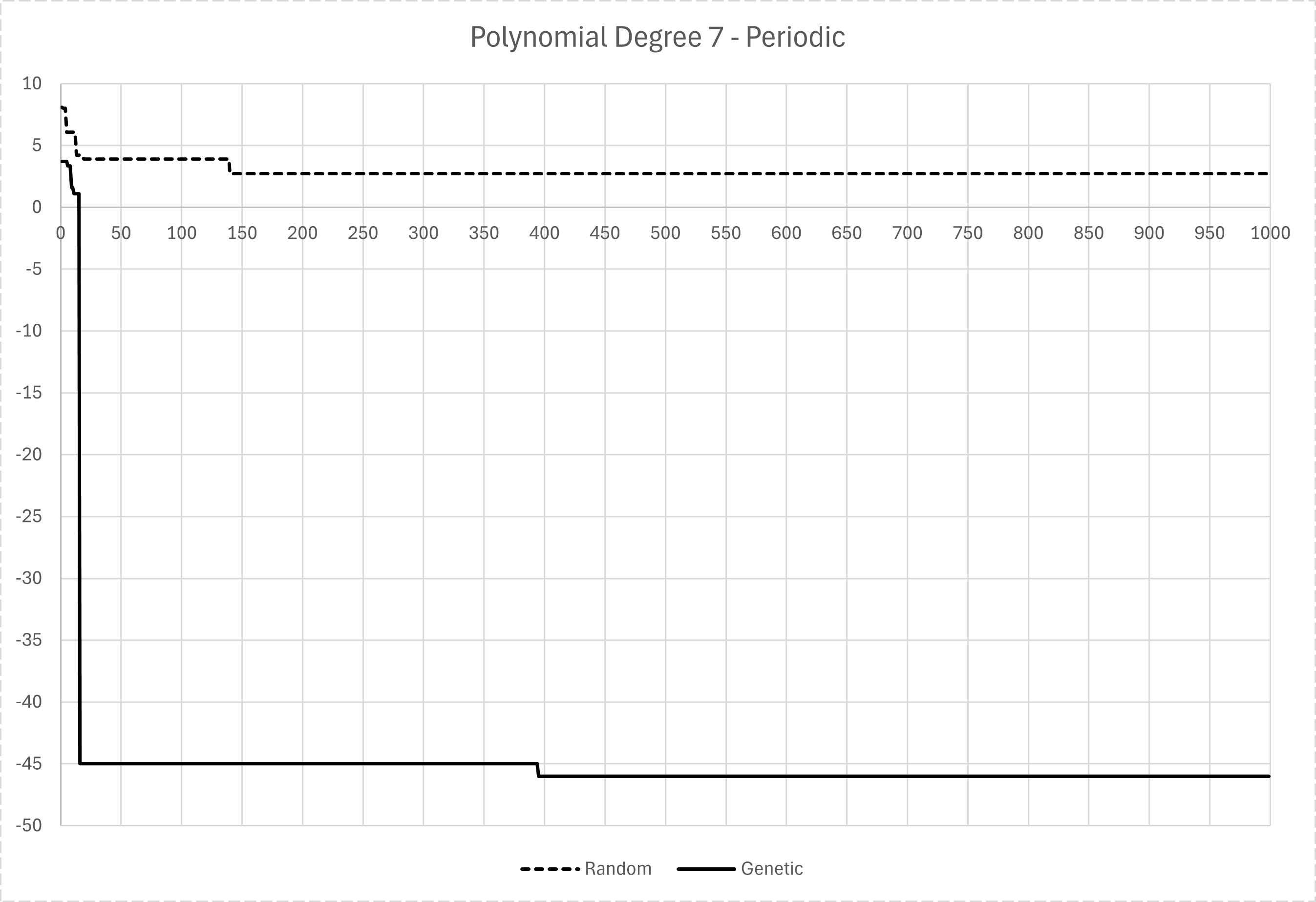}
\end{center}

\begin{center}
\includegraphics[scale=0.45]{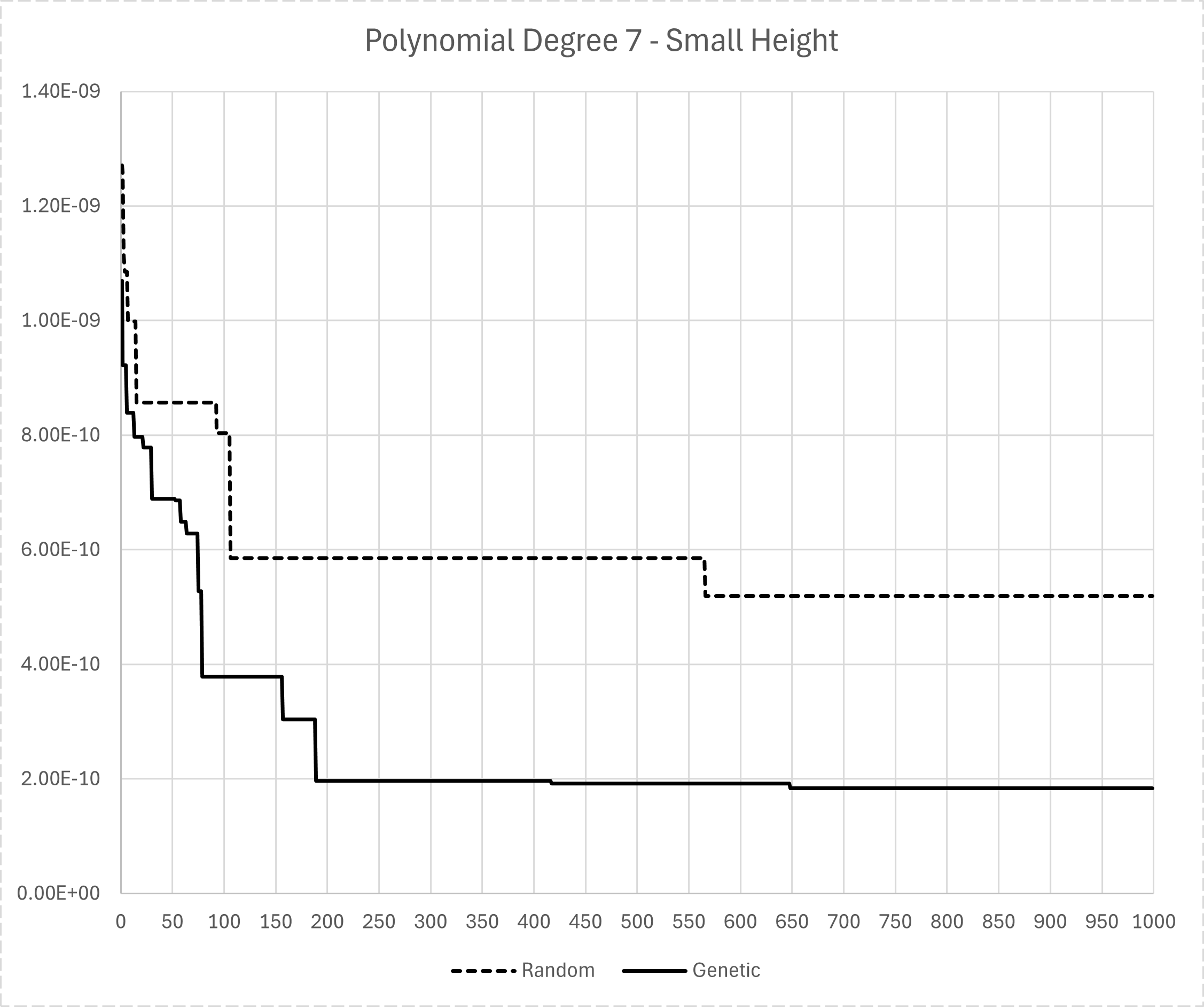}
\end{center}

There are a few clear next steps from this work. We enumerate just a few broad questions.
\begin{enumerate}
    \item Can the fitness function be modified to improve performance on the Many Preperiodic Points problem (Morton-Silverman Uniform Boundedness)?
    \item Can the output of this algorithm be used as ``seed'' values for a more advanced machine learning algorithm to improve the extreme examples found?
    \item Can a more advanced version of machine learning push past the degree limitations of this genetic algorithm?
    \item Can the algorithm be modified to work over finite fields and examine extreme properties of the associated preperiodic point graph structure?
\end{enumerate}

\subsection{Interesting Findings}
A final step we performed was to look at the final data sets through the human perspective and see if there were any patterns or families that could be exploited for general results. We mainly identified only a single reoccuring pattern. Doyle and Hyde \cite{Doyle-Hyde} introduced a notion of \emph{dynamical compression} in their search for sets of maps that shared many common preperiodic points. The idea is that a map $f$ exhibits dynamical compression if it maps a consecutive interval of integers onto itself, i.e., if there exists an $m$ so that
\begin{equation*}
    f(\{1,2,\ldots,m\}) \subseteq \{1,2,\ldots,m\}.
\end{equation*}
With an appropriate conjugation, the classic example $f(z) = z^2 - \frac{29}{16}$ exhibits this dynamical compression \cite[Example 1.6]{Doyle-Hyde}. They also give several cubic polynomials that match the best known results of Benedetto et al \cite{Benedetto4} that exhibit dynamical compression \cite[Example 1.7]{Doyle-Hyde}. They further study a family of maps with dynamical compression to prove their lower bounds on the number of rational preperiodic points. 

In our results, we found an extraordinarily high number of extreme examples exhibiting this dynamical compression (i.e., the orbit of $0$ consisted of a permutation of a consecutive sequence of integers). It is not clear why these types of maps are exhibiting such extreme behavior nor the frequency with which they are occurring. Furthermore, we did not investigate whether the maps that were not obviously exhibiting this compression could be conjugated to a form that did (as with the $z^2 - \frac{29}{16}$ example). It would be interesting to further explore this topic.

\section{Data - Examples} \label{sect_data}
\def\arraystretch{1.5}%

The following tables summarize the best results from running this algorithm on the four specific problems.

\subsection{Small Height Ratio}
For quadratic and cubic polynomials and quadratic rational functions our algorithm matched the best known values.

By choosing an orbit of length $d+1$ (or $2d+2$) which contains points of small height, we are already achieving a much lower height ratio than a generic orbit that should grow approximately by a factor of $d$ with each iterate. There is no obviously determinable pattern from our data, but we give an approximate growth factor of the exponent of the ratio in terms of $d$.

\begin{enumerate}
    \item The smallest height ratio seems to be growing approximately as $10^{-\frac{3d}{2}}$ for polynomial maps.
    \item The smallest height ratio seems to be growing approximately as $10^{-2d}$ for rational maps.
\end{enumerate}

{\tiny
\begin{center}
\begin{tabular}{|c|c|c|c|c|c|}
\hline
\multicolumn{6}{|c|}{Small Height ratio}\\
\hline
Degree & type & Orbit & result & Best Known & citation\\
\hline
2 & Polynomial & [0, -2,1,-3] & $0.006604$ & $0.006604$ & \cite{Benedetto8}\\
\hline
3 & Polynomial & $[0, 1, -3, -4, -8]$ & $0.000092099$ & $0.000092099$ & \cite{Benedetto4}\\
\hline
4 & Polynomial & $[0, 1, 4, 5, -1, 3]$ & $2.9015 \cdot 10^{-6}$ &  & \\
\hline
5 & Polynomial & $[0, -1, -2, -6, -4, -5, -3]$ & $9.1519\cdot 10^{-9}$ &  & \\
\hline
6 & Polynomial & $[0, 1, -1, -4, -7, -3, -6, -9]$ & $1.8372 \cdot 10^{-9}$ &  & \\
\hline
7 & Polynomial & $[0, -1, 1, 11, 2, 10, 8, 9, 3]$ & $1.0564 \cdot 10^{-10}$ &  & \\
\hline
8 & Polynomial & $[0, 1, 3, -10, -12, 2, -1, -9, -11, -8]$ & $8.7710 \cdot 10^{-13}$ &  & \\
\hline
9 & Polynomial & $[0, -1, -12, -9, 1, -4, -3, -8, -11, -10, 2]$ & $1.7173\cdot 10^{-14}$ &  & \\
\hline
10 & Polynomial & $[0, [-1, 3, 15, 9, 12, 1, 10, 11, 4, 14, 2]$ & $9.9368\cdot 10^{-16}$ &  & \\
\hline
11 & Polynomial & $[0, -1, 1, -10, -9, -12, -11, -8,$ & $7.3446\cdot 10^{-17}$ &  & \\
&&$ -14, -3, -2, 2, -13]$ &&&\\
\hline
12 & Polynomial & $[0, 1, -2, -11, -8, -1, 12, 7, 3,$ & $5.9715 \cdot 10^{-19}$ &  & \\
&&$ -18, 8, 5, 2, -15]$ &&&\\
\hline
\hline
2 & Rational & $[0, -1, -16, 4, 8, 2]$& $0.0004657$ & $$0.0004657$$ & \cite{Benedetto8}\\
\hline
3 & Rational & $[0, -1, -7, 3, -11, 5, -3, -5]$& $3.0790 \cdot 10^{-6}$ &  & \\
\hline
4 & Rational & $[0, -1, -4, -3, -9, 2, -6, -2, -5, 1]$ & $2.1843\cdot 10^{-8}$ &  & \\
\hline
5 & Rational & $[0, 1, -1, 7, -17, 13, -27, -8, -2, 14, -29, 4]$ & $3.6941 \cdot 10^{-10}$ &  & \\
\hline
\end{tabular}
\end{center}
}

\subsection{Many Preperiodic Points (Morton-Silverman Uniform Boundedness)}
For low degree functions, our algorithm is close to the best known bound, but appears to be getting further from the optimal values as the degree increases. As discussed above, this appears to mainly be due to the algorithm converging to long sparse cycles rather than dense preperiodic graphs as in the best known examples. Note that the counts listed here include the fixed point at infinity.

{\tiny
\begin{center}
\begin{tabular}{|c|c|c|c|c|c|}
\hline
\multicolumn{6}{|c|}{Many Preperiodic Points}\\
\hline
Degree & type & Orbit & result & Best Known & citation\\
\hline
2 & Polynomial & $[0,1,-1,2]$ & 9 & $9$ & \cite{Poonen}\\
\hline
3 & Polynomial & $[0, 1, -3, -1, 4]]$ & $10$ & $11$ & \cite{Benedetto4}\\
\hline
4 & Polynomial & $[0, 1, 10, 4, 8, 3]$ & $11$ & $11$ & \cite{Doyle-Hyde}\\
\hline
5 & Polynomial & $[0, 1, 3, 5, -7, -9, -11]$ & $12$ & $14$ & \cite{Doyle-Hyde}\\
\hline
6 & Polynomial & $[0, 1, -9, -2, -8, -3, -7, -1]$ & $12$ & $15$ & \cite{Doyle-Hyde}\\
\hline
7 & Polynomial & $[0, 1, 2, -8, -1, -4, -6, -2, -7]$ & $13$ & $16$ & \cite{Doyle-Hyde}\\
\hline
8 & Polynomial & $[0, -1, 3, 2, -2, 1, -3, -6, -4, -7]$ & $12$ & $17$ & \cite{Doyle-Hyde}\\
\hline
\hline
2 & Rational & $[0, -1, -3, -6, -2, -4]$ & $14$ & $14$ & \cite{Benedetto8}\\
\hline
3 & Rational & $[0, 1, 2, 3, 5, 4, -7, 6]$ & $13$ &  & \\
\hline
4 & Rational & $[0, 1, 3, 2, 5, -2, -3, -4, -5, -9]$ & $11$ &  & \\
\hline
\end{tabular}
\end{center}
}

\subsection{Long Periodic Cycles}
For quadratic and cubic polynomials, our algorithms attains the best known values. For quadratic rational functions we find a period 6 $\Q$-rational periodic point, whereas the best known is period $7$.

It is interesting to note that the longest known periodic cycle increases linearly with the degree. In particular, we make the following formalizations.
\begin{itemize}
    \item For $f$ a degree $d\geq 2$ polynomial map, there were no examples of a $\Q$-rational periodic point with minimal period larger than $d+3$.
    \item For $f$ a degree $d\geq 2$ rational map, there were no examples of a $\Q$-rational periodic point with minimal period larger than $2d+3$.
\end{itemize}

{\tiny
\begin{center}
\begin{tabular}{|c|c|c|c|c|c|}
\hline
\multicolumn{6}{|c|}{Long Periodic Cycles}\\
\hline
Degree & type & Orbit & result & Best Known & citation\\
\hline
2 & Polynomial & $[0,1,-1,2]$ & 3 & $3$ & \cite{Poonen}\\
\hline
3 & Polynomial & $[0, 1, -1, 3, -3]$ & $5$ & $5$ & \cite{Benedetto4}\\
\hline
4 & Polynomial & $[0,1, 5, 4, 6, 2]$ & $6$ &  & \\
\hline
5 & Polynomial & $[0,1, 7, 5, -1, 3, 8]$ & $8$ &  & \\
\hline
6 & Polynomial & $[0,1, -8, -5, -1, 2, -7, -4]$ & $8$ &  & \\
\hline
7 & Polynomial & $[0,-1, 1, -2, -6, -8, -5, -4, -3]$ & $10$ &  & \\
\hline
8 & Polynomial & $[0,1, 7, 5, 8, 6, 9, 4, 2, 3]$ & $11$ &  & \\
\hline
9 & Polynomial & $[0,1, 3, 2, 4, 8, 9, 10, 6, 7, 5]$ & $11$ &  & \\
\hline
10 & Polynomial & $[0,1, 2, 3, -6, -11, -2, -7, 4, -9, -10, -1]$ & $12$ &  & \\
\hline
11 & Polynomial & $[0,1, 7, 4, -2, 11, 6, -4, -1, 12, 9, 8, 10]$ & $13$ &  & \\
\hline
12 & Polynomial & $[0,-1, -6, 7, -7, -8, 4, 6, -3, 5, -4, 1, 3, -5]$ & $14$ &  & \\
\hline
13 & Polynomial & $[0,-1, -3, -2, -13, 1, -12, -10,$ & $15$ &  & \\
 & & $-11, -5, -7, -8, -6, -9, -4]$ &&&\\
\hline
\hline
2 & Rational & $[0, 1, -1, -9, -14, 26]$ & $6$ & $7$ &  \cite{Benedetto8}\\
\hline
3 & Rational & $[0,1, 5, -9, -5, -3, -1, 3]$ & $8$ &  & \\
\hline
4 & Rational & $[0,-1, -2, -4, 3, 6, 5, 4, -5, 1]$ & $10$ &  & \\
\hline
\end{tabular}
\end{center}
}

\subsection{Long Preperiodic Tails}
Our algorithm attains the best known values for all previously explored cases.
It is interesting to note that the longest known rational tail increases linearly with the degree. In particular, we make the following formalizations.
\begin{itemize}
    \item For $f$ a degree $d\geq 2$ polynomial map there were no examples of a $\Q$-rational preperiodic point with minimal tail larger than $d+2$. 
    \item For $f$ a degree $d\geq 2$ rational map there were no examples of a $\Q$-rational preperiodic point with minimal tail larger than $2d+2$.
\end{itemize}

{\tiny
\begin{center}
\begin{tabular}{|c|c|c|c|c|c|}
\hline
\multicolumn{6}{|c|}{Long Preperiodic Tail}\\
\hline
Degree & type & Orbit & result & Best Known & citation\\
\hline
2 & Polynomial & $[0,1,-1,2]$ & 2 & $2$ & \cite{Poonen}\\
\hline
3 & Polynomial & $[0, 1, -1, 2, -3]$ & $4$ & $4$ & \cite{Benedetto4}\\
\hline
4 & Polynomial & $[0,1, -1, 2, -3, 4]$ & $5$ &  & \\
\hline
5 & Polynomial & $[0,-1, -2, 3, 5, 6, 7]$ & $6$ &  & \\
\hline
6 & Polynomial & $[0,1, 2, 7, 6, 5, 8, 9]$ & $7$ &  & \\
\hline
7 & Polynomial & $[0,1, 3, 5, 4, 6, -1, -3, -4]$ & $9$ &  & \\
\hline
8 & Polynomial & $[0,-1, -5, -7, -3, -8, -2, -6, -4, -9]$ & $9$ &  & \\
\hline
9 & Polynomial & $[0,1, -1, 6, 11, 8, 9, 10, 4, 3, 2]$ & $10$ &  & \\
\hline
10 & Polynomial & $[0,-1, 2, 4, -7, -5, 1, -4, -3, -8, 3, -6]$ & $11$ &  & \\
\hline
11 & Polynomial & $[0,-1, -5, -7, -2, -8, -11, 3, -6, 1, -10, -9, 2]$ & $12$ &  & \\
\hline
12 & Polynomial & $[0,1, 10, 13, 3, 6, 9, 2, 5, 12, -1, 4, 11, -2]$ & $13$ &  & \\
\hline
13 & Polynomial & $[0,-1, -2, -4, -12, -14, -8,-10,$ & $14$ &  & \\
&&$-3, -5, -7, -13, -11, 1, -6]$&&&\\
\hline
\hline
2 & Rational & $[0,-1, -3, -6, -12, 3]$ & $6$ & $6$ & \cite{Benedetto8}\\
\hline
3 & Rational & $[0,-1, 2, 3, -3, 11, 1, 5]$ & $8$ &  & \\
\hline
4 & Rational & $[0,-1, 1, 9, 5, -11, -9, -7, -5, -3]$ & $9$ &  & \\
\hline
\end{tabular}
\end{center}
}

%\bibliography{biblio}
%\bibliographystyle{plain}
\providecommand\biburl[1]{\texttt{#1}}

\begin{all_data}
\newpage
\section{Appendix - Extended Data Set}

The data is arranged first by degree, then by problem. Each table contains the following information.
\begin{enumerate}
    \item The orbit which determines the map

    \item The value obtained by this map for this problem.
    
    \item An explicit form for the map from interpolating the orbit.
    
    \item A reduced representation in the conjugacy class in the terminology of Hutz-Stoll \cite{Hutz20}. This is a map with minimal resultant and smallest coefficients (in terms of height). This model is omitted when it was not reasonably computable.
    
    \item The conjugation from the interpolated model to the reduced model.
\end{enumerate}

%%%%%%%%%%%%%%%%%%%%%%%%%%%%%%%%%%%%%
\subsection{Degree 2: Polynomial}

{\small
\begin{center}
\begin{tabular}{|l|c|c|c|c|}
\hline
Orbit & result & Interpolated Map & Normal Form & conjugation \\
\hline
\hline
\multicolumn{5}{|c|}{Small Height ratio}\\
\hline
[0,-2,1,-3] & $0.006604$ & $1/6z^2 - 7/6z - 2$ &  $z^2 - \frac{181}{144}$ &
$\sm{6,7/2,0,1}$  \\
\hline
[0,1,-1,2] & $0.01102$ & $-1/6z^2 - 17/6z + 1$ &  $z^2 - \frac{517}{144}$ &
$\sm{-6,-17/2,0,1}$  \\
\hline
[0,18,21,14] & $0.01346$ & $-5/42z^2 + 97/42z + 18$ &  $z^2 - \frac{16381}{7056}$ &
$\sm{-42/5,97/10,0,1}$  \\
\hline
\hline
\multicolumn{5}{|c|}{Many Preperiodic Points}\\
\hline
[0,1,-1,2] & 9 (2,3)& $-1/2z^2 - 3/2z + 1$ &  $z^2 - \frac{29}{16}$ &
$\sm{-2,-3/2,0,1}$  \\
\hline
[0,1,4,7] & 9 (1,3) & $-1/2z^2+7/2z +1$ & $z^2-\frac{29}{16}$ & $\sm{-2,7/2,0,1}$  \\
\hline
[0,-4,2,-1] & 9 (1,3) & $2/3x^2 + 7/6x - 4$ & $x^2 - \frac{439}{144}$ & $\sm{3/2,-7/8,0,1}$  \\
\hline
\hline
\multicolumn{5}{|c|}{Long Periodic Cycles}\\
\hline
[0,1,-1,2] & (2,3) & $-1/2z^2 - 3/2z + 1$ &  $z^2 - \frac{29}{16}$ &
$\sm{-2,-3/2,0,1}$ \\
\hline
[0,1,4,7] & (1,3) & $-1/2z^2+7/2z +1$ & $z^2-\frac{29}{16}$ &$\sm{-2,7/2,0,1}$ \\
\hline
\hline
\multicolumn{5}{|c|}{Long Preperiodic Tail}\\
\hline
[0,1,-1,2] & (2,3)& $-1/2x^2 - 3/2x + 1$ & $x^2 - \frac{29}{16}$ & $\sm{-2,-3/2,0,1}$  \\
\hline
[0,4,2,1] & (2,2) & $1/2x^2 - 5/2x + 4$ & $x^2 - \frac{13}{16}$ & $\sm{2,5/2,0,1}$  \\
\hline
[0,1,-1,4]& (2,2) & $1/2x^2 - 5/2x + 1$ & $x^2 - \frac{37}{16}$ & $\sm{2,5/2,0,1}$ \\
\hline
[0,1,-3,16] & (2,2) & $1/4x^2 - 17/4x + 1)$ & $x^2 - \frac{409}{64}$ & $\sm{4,17/2,0,1}$  \\
\hline
\end{tabular}
\end{center}
}

\begin{comment}
Code to analyze maps....
\begin{verbatim}
F=orbit_to_map_polynomial([3,1,5],2)
print(F.dehomogenize(1))
G,m,phi=F.normal_form(return_conjugation=True)
G.dehomogenize(1), m

h_P=F.canonical_height(F.domain()(0),error_bound=0.000001)
h_F = max([sig.global_height() for sig in F.sigma_invariants(1)])
h_P/h_F

h_Q=G.canonical_height(G.domain()(7/12),error_bound=0.000001)
h_G = (181/144).global_height()
h_Q/h_G
\end{verbatim}

\subsubsection{Best Previously Known}

\begin{center}
\begin{tabular}{|l|c|c|}
\hline
Map & value & citation\\
\hline
\multicolumn{3}{|c|}{Long Periodic Cycles}\\
\hline
$z^2-29/16$& (2,3)& Poonen\\
\hline
\multicolumn{3}{|c|}{Long Preperiodic Tail}\\
\hline
$z^2-29/16$& (2,3)& Poonen\\
\hline
\multicolumn{3}{|c|}{Many Preperiodic Points}\\
\hline
$z^2-29/16$& 9 (2,3) & Poonen\\
\hline
\multicolumn{3}{|c|}{Small Height ratio}\\
\hline
 $z^2 - \frac{181}{144}$ & $0.006604$ & Benedetto et al\\
\hline
\end{tabular}
\end{center}

\end{comment}

\subsection{Degree 2: Rational Functions}
{\small
\begin{center}
\begin{tabular}{|l|c|c||c|c|}
\hline
Orbit & result & Interpolated Map & Reduced Form & conjugation  \\
\hline
\multicolumn{5}{|c|}{Small Height ratio}\\
\hline
$[0,-1, -16, 4, 8, 2]$& $0.0004657$ & $\frac{592z^2 - 3424z + 1024}{173z^2 - 536z - 1024}$ & $\frac{7z^2 - z - 6}{-z^2 + 9z + 6}$ & $\sm{-8,16,-7,-1}$ \\
%&&& Rob: $\frac{10z^2 - 7z - 3}{10z^2 + 37z + 9}$ & $\sm{0,-16,15,1}$ & $z= \infty$ & \\
\hline
%$[1, 3, 9, 5, -15]$&same as the first & & & \\
%\hline
$[0,1, 6, -6, -5, -3]$ & $0.0011830$ & $\frac{285z^2 + 2085z + 990}{-79z^2 - 351z + 990}$ & $\frac{75z^2 - 475z - 440}{-9z^2 + 134z - 440}$ & $\sm{-3,0,1,-4}$  \\
%&&& Rob: $\frac{1375z^2 - 847z - 528}{1375z^2 + 4025z + 1320}$ & $\sm{0,6,5,1}$ & $z= \infty$ & \\
\hline
\end{tabular}

\begin{tabular}{|l|c|c|c||c|}
\hline
Orbit & result & Interpolated Map & Reduced Form & conjugation \\
\hline
\multicolumn{5}{|c|}{Many Preperiodic Points}\\
\hline
$[0,-1, -3, -6, -2, -4] $  & 14,(6, 2) & $\frac{7z^2 + 31z - 6}{-2z^2 - 6z + 6)}$& $\frac{-15z^2 + 20z - 65}{42z^2 - 3z - 75)}$ & $\sm{-3,2,0,-1}$  \\
%&&& Rob: $\frac{8z^2 - 5z - 3}{8z^2 + 70z + 12}$ & $\sm{0,-3,2,1}$ & $z= -1$ & \\
\hline
$[0,-1, 2, 6, 4, 5]$ & 14,(1, 5)& $\frac{83z^2 - 403z - 66}{22z^2 - 122z + 66}$& $\frac{39z^2 - z - 40}{66z^2 - 34z - 30}$ & $\sm{-3,-2,0,-1}$  \\
%&&& Rob: $\frac{99z^2 - 119z + 20}{99z^2 - 469z + 90}$ & $\sm{0,2,-3,1}$ & $z= \infty, \frac{5}{9}$ & \\
\hline
$[0,1, -1, 3, -6, -15]$ & 14,(3, 3) & $\frac{9z^2 + 84z - 45}{z^2 - 4z - 45}$ &  $\frac{5z - 5}{z^2 - 2z - 5}$ & $\sm{-3,-9,1,-1}$ \\
%&&& Rob: $\frac{45z^2 - 65z + 20}{45z^2 - 87z + 30}$ & $\sm{0,-1,-2,1}$ & $z= \infty, \frac{1}{3}$ & \\
\hline
$[0,1, -1, 4, -5, -14]$ & 14,(3, 3) &  $\frac{3z^2 + 29z - 14}{-4z - 14}$& $\frac{3z^2 - 7z - 6}{-4z + 2}$ & $\sm{3,-5,0,1}$  \\
%& && Rob: $\frac{56z^2 - 81z + 25}{56z^2 - 114z + 40}$ & $\sm{0,-1,-2,1}$ & $z= \infty, \frac{5}{14}$ & \\
\hline
\end{tabular}

\begin{tabular}{|l|c|c|c|c|}
\hline
Orbit & result & Interpolated Map & Reduced Form & conjugation \\
\hline
\hline
\multicolumn{5}{|c|}{Long Periodic Cycles}\\
\hline
\hline
$[0,1, -1, -9, -14, 26]$ & (0, 6) &  $\frac{83z^2 - 2179z + 546}{73z^2 + 931z + 546}$ & $\frac{-195z^2 - z + 124}{-195z^2 + 117z}$ & $\sm{5,4,-5,4}$  \\
%&&& Rob: $\frac{1092z^2 - 1716z + 624}{1092z^2 - 3271z + 1404}$ & $\sm{0,-1,-2,1}$ &  $z= \infty$& \\
\hline
$[0,-1, -4, -10, -2, 2]$ & (1, 5)  & $\frac{28z^2 + 72z + 80}{-19z^2 - 90z - 80}$ & $\frac{3z^2 - 13z + 14}{3z^2 + 3z}$ & $\sm{-2,-2,2,-1}$ \\
%&&& Rob: $\frac{15z^2 - 14z - 1}{15z^2 - 8z + 5}$ & $\sm{0,-4,3,1}$ & $z= \infty$& \\
\hline
$[0,-1, 1, 9, -3, 4]$ & (1, 5)  & $\frac{91z^2 - 208z - 99}{-11z^2 - 112z + 99}$ & $\frac{10z^2 + 45z - 55}{-2z^2 - 24z - 40}$ & $\sm{-1,-10,1,0}$  \\
%&&& Rob: $\frac{99z^2 - 123z + 24}{99z^2 - 203z + 54}$ & $\sm{0,1,-2,1}$ & $z= \infty$& \\
\hline
$[0,1, 3, 5, 15, -3]$ & (1, 5)  & $\frac{27z^2 + 354z - 45}{-29z^2 + 186z - 45}$ &  $\frac{21z + 21}{-12z^2 + 19z + 21}$ & $\sm{-3,-3,1,-3}$ \\
%&&& Rob: $\frac{5z^2 - 12z + 7}{5z^2 - 54z - 35}$ & $\sm{0,3,2,1}$ & $z= \infty$& \\
\hline
\end{tabular}

\begin{tabular}{|l|c|c|c||c|}
\hline
Orbit & result & Interpolated Map & Reduced Form & conjugation \\
\hline
\hline
\multicolumn{5}{|c|}{Long Preperiodic Tail}\\
\hline
$[0,-1, -3, -6, -12, 3]$ & (6, 2)  & $\frac{129z^2 + 1737z - 702}{-19z^2 - 87z + 702}$ &  $\frac{44z^2 - 9z - 38}{22z^2 + 58z - 24}$ & $\sm{-3,-6,1,0}$\\
%&&& Rob: $\frac{52z^2 - 30z - 22}{52z^2 + 245z + 88}$ & $\sm{0,-3,2,1}$ & $z= \infty$& \\
\hline
$[0,-1, -9, 27, 3, -15]$ & (6, 2)  & $\frac{63z^2 - 11934z + 1215}{157z^2 + 510z - 1215}$ & $\frac{7z - 7}{-5z^2 + 4z + 7}$ & $\sm{43,-63,11,7}$ \\
%&&& Rob: $\frac{10z^2 - 6z - 4}{10z^2 + 113z + 24}$ & $\sm{0,-9,8,1}$ & $z= \infty$ & \\
\hline
$[0,1, 3, -8, -7, -11]$ &(6, 2)  & $\frac{43z^2 + 106z + 175}{-14z^2 - 53z + 175}$ &  $\frac{15z^2 + 20z - 65}{42z^2 - 3z - 75}$ & $\sm{-3,2,0,-1}$\\
%&&& Rob: $\frac{700z^2 - 95z - 605}{700z^2 + 1336z + 880}$ & $\sm{0,3,2,1}$ & $z= \infty$ & \\
\hline
\end{tabular}

\end{center}
}

\subsection{Degree 3: Polynomial}
{\small
\begin{center}
\begin{tabular}{|l|c|c|c|}
\hline
Orbit & result & Interpolated Map & conjugation \\
\hline
\hline
\multicolumn{4}{|c|}{Small Height ratio}\\
\hline
$[0,1, -3, -4, -8]$ & 0.000092099 &$-1/6z^3 - 7/4z^2 - 25/12z + 1$\\
\cline{2-4}
& normal form& $-25/24z^3 + 97/24z + 1$ & $\sm{-2,-7,0,5}$ \\
\cline{2-4}
& reduced form & $-1/6z^3 - 7/4z^2 - 25/12z + 1$ & $\sm{1,0,0,1}$ \\
\hline
$[0,1, 3, -1, -3]$&0.00016547 & $-1/12z^3 - z^2 + 37/12z + 1$&  \\
\cline{2-4}
& normal form& $-27z^3 + 85/12z + 1$  & $\sm{-1,-4,0,2}$ \\
\cline{2-4}
& reduced form & $-1/3x^3 - 5/2x^2 + 5/6x + 1$ & $\sm{2,1,0,1}$ \\
\hline
$[0,1, 2, 4, 3]$&0.00016738& $-1/3z^3 + 3/2z^2 - 1/6z + 1$  & \\
\cline{2-4}
& normal form& $-3/4z^3 + 25/12z + 1$ & $\sm{2,-3,0,3}$ \\
\cline{2-4}
& reduced form & $-1/3z^3 + 3/2z^2 - 1/6z + 1$ & $\sm{1,0,0,1}$ \\
\hline
\hline
\multicolumn{4}{|c|}{Many Preperiodic Points}\\
\hline
$[0,1, -3, -1, 4]$ & 10 (3,2) & $1/3z^3 - 1/2z^2 - 23/6z + 1$ & \\
\cline{2-4}
& normal form& $3/4z^3 - 49/12z + 1$  & $\sm{-2,1,0,3}$ \\
\cline{2-4}
& reduced form & $1/3z^3 - 1/2z^2 - 23/6z + 1$ & $\sm{1,0,0,1}$\\
\hline
\hline
\multicolumn{4}{|c|}{Long Periodic Cycles}\\
\hline
\hline
$[0,1, -1, 3, -3]$ & (1, 5) & $-1/6z^3 - 7/4z^2 - 25/12z + 1$ & \\
\cline{2-4}
 & normal Form& $-25/24z^3 + 97/24z + 1$ & $\sm{-2,-7,0,5}$ \\
 \cline{2-4}
 & reduced form & $1/3x^3 + 1/2x^2 - 11/6x - 1$ & $\sm{2, 1, 0, 1}$\\
\hline
$[0,1, -1, 2, -2]$ & (0, 5) & $1/3z^3 - 1/2z^2 - 11/6z + 1$ & \\\cline{2-4}
& normal form&$1/12z^3 - 25/12z + 1$ & $\sm{-2,1,0,1}$\\
\cline{2-4}
& reduced form & $1/3z^3 - 1/2z^2 - 11/6z + 1$ & $\sm{1,0,0,1}$\\
\hline
\hline
\multicolumn{4}{|c|}{Long Preperiodic Tail}\\
\hline
$[0,1, -1, 2, -3]$ & (4, 1) & $1/6z^3 - 1/2z^2 - 5/3z + 1$  & \\\cline{2-4}
& normal form& $2/3z^3 - 13/6z + 1$ & $\sm{-1,1,0,2}$\\
\cline{2-4}
 & reduced form & $1/6z^3 - 1/2z^2 - 5/3z + 1$ & $\sm{1,0,0,1}$\\
\hline
$[0,1, 4, 9, 10]$& (4, 1) & $1/60z^3 - 5/12z^2 + 17/5z + 1$  & \\
\cline{2-4}
& normal form& $76729/1574640z^3 - 13/180z + 1$ & $\sm{162, -1350, 0, 277}$ \\
\cline{2-4}
 & reduced form & $4/15x^3 + 1/3x^2 + 1/15x$ & $\sm{4,10,0,1}$ \\
\hline
$[0,1, 7, 3, 8]$& (4, 1) & $37/168z^3 - 19/7z^2 + 1427/168z + 1$ & \\
\cline{2-4}
& normal Form& $11916147/34747958z^3 - 2359/888z + 1$  & $\sm{9583, -39368,0,11958}$ \\
\cline{2-4}
& reduced form & $37/42x^3 + 33/28x^2 - 179/84x - 13/14$ & $\sm{2,5,0,1}$ \\
\hline
\end{tabular}
\end{center}
}

\begin{comment}
Code to analyze maps....
\begin{verbatim}
R.<a,b>=QQ[]
P.<x,y>=ProjectiveSpace(R,1)
F=DynamicalSystem([a*x^3 + b*x*y^2, y^3])
sig1 = F.sigma_invariants(1)
G=DynamicalSystem([a*x^3 + b*x*y^2+y^3, y^3])
sig2 = G.sigma_invariants(1)

F=orbit_to_map_polynomial([1,-3,-4,-8],3)
print(F.dehomogenize(1))
P.<x,y>=ProjectiveSpace(QQ,1)
sig = F.sigma_invariants(1)
I = R.ideal([R(sig[i] - sig1[i]) for i in range(len(sig))])
V1= I.variety()
print(V1)
I = R.ideal([R(sig[i] - sig2[i]) for i in range(len(sig))])
V2=I.variety()
V2

b1 = V2[0][b]
a1 = V2[0][a]
A.<z>=AffineSpace(QQ,1)
f=DynamicalSystem_affine([a1*z^3+b1*z+1]).homogenize(1)
f.change_ring(QQbar).conjugating_set(F.change_ring(QQbar))


hF = max([t.global_height() for t in f[0].coefficients()])
hs = max([t.global_height() for t in F.sigma_invariants(1)])
hQ=F.canonical_height(P(0), error_bound=0.0000001)
hQ, hF, hs, hQ/hF, hQ/hs
\end{verbatim}

\subsubsection{Best Previously Known}
\begin{center}
\begin{tabular}{|l|c|c|}
\hline
Map & value & citation\\
\hline
\multicolumn{3}{|c|}{Long Periodic Cycles}\\
\hline
$1/12z^3 - 25/12z + 1$; $z=0$& $(1,5)$ &Bendetto (2007)\\
\hline
\multicolumn{3}{|c|}{Long Preperiodic Tail}\\
\hline
$-4/3z^3 + 37/12z + 1$; $z=-5/4$ & (4,1) 7 & Bendetto(2007)\\
\hline
\multicolumn{3}{|c|}{Many Preperiodic Points}\\
\hline
$-3/2z^3 + 19/6z$ & 11 (2,2)& Benedetto (2007)\\
\hline
\multicolumn{3}{|c|}{Small Height ratio}\\
\hline
$-25/24z^3+97/24z + 1$; $z=-7/5$ & 0.0002513 & Benedetto (2007)\\
 & sigma: $0.000092089$ & \\
\hline
\end{tabular}
\end{center}
\end{comment}

\subsection{Degree 3: Rational Functions}
{\small
\begin{center}

\begin{tabular}{|l|c|c|c|}
\hline
Orbit & result & Interpolated Map  & conjugation \\
\hline
\multicolumn{4}{|c|}{Small Height ratio}\\
\hline
$[0,-1, -7, 3, -11, 5, -3, -5]$& $3.0790 \cdot 10^{-6}$ & $\frac{1115z^3 - 2405z^2 - 37855z + 3465}{201z^3 + 173z^2 + 1907z - 3465}$ &  \\
\cline{2-4}
&reduced & $\frac{658z^3 - 1545z^2 - 3442z + 2025}{201z^3 - 215z^2 + 541z - 675}$ & $\sm{2,-1,0,1}$   \\
\hline
$[0,-1, 1, -25, 3, 7, 11, 15]$ & $3.25839\cdot 10^{-6}$ & $\frac{183z^3 - 3203z^2 + 15445z + 5775}{57z^3 - 1117z^2 + 6107z - 5775}$ & \\
\cline{2-4}
&reduced & $\frac{28z^3 + 72z^2 + 20z - 57}{-28z^3 + 87z^2 - 5z - 54}$ & $\sm{7,-4,1,0}$  \\
\hline
$[0,1, -6, -1, -4, -3, 6, 3]$ & $3.29835\cdot 10^{-6}$& $\frac{21z^3 + 312z^2 + 579z - 72}{26z^3 + 47z^2 - 141z - 72}$ &\\
\cline{2-4}
&reduced & $\frac{32z^3 + 148z^2 - 253z - 72}{-8z^3 - 217z^2 - 98z + 48}$ & $\sm{0,-3,1,1}$  \\
\hline
$[0,-1, -15, -3, -5, -7, 5, 1]$& $3.5587 \cdot 10^{-6}$ & $\frac{50z^3 + 195z^2 - 1640z - 525}{-12z^3 + 9z^2 + 630z + 525}$ &\\
\cline{2-4}
&reduced & $\frac{60z^3 + 63z^2 - 117z - 48}{30z^3 + 44z^2 + 24z - 14}$ & $\sm{-3,-4,1,0}$  \\
\hline
$[0,1, 15, 3, 5, 7, -5, -1]$ & $3.55891 \cdot 10^{-6}$ & $\frac{50z^3 - 195z^2 - 1640z + 525}{12z^3 + 9z^2 - 630z + 525}$ &\\
\cline{2-4}
& reduced & $\frac{60z^3 - 63z^2 - 117z + 48}{-30z^3 + 44z^2 - 24z - 14}$ & $\sm{3,-4,1,0}$   \\
\hline
$[0,1, -7, 7, -3, -5, 5, -1]$& $3.68152 \cdot 10^{-6}$ & $\frac{153z^3 - 185z^2 - 1977z + 665}{x^3 - 185z^2 - 289z + 665}$ &\\
\cline{2-4}
& reduced & $\frac{76z^3 + 114z^2 - 154z - 96}{x^3 - 91z^2 - 164z + 24}$ & $\sm{2,1,0,1}$   \\
\hline
$[0,-1, 1, -3, 2, 3, 5, 7]$&$ 3.70697\cdot 10^{-6}$ & $\frac{6z^3 - 125z^2 + 126z + 65}{4z^3 - 43z^2 + 80z - 65}$ &\\
\cline{2-4}
& reduced & $\frac{-2z^3 + 38z^2 + 28z - 12}{-8z^3 + 31z^2 - 3z + 6}$ & $\sm{2,1,0,1}$   \\
\hline
$[0,-1, -5, -11, 4, 7, 9, 10]$ & $3.91358 \cdot 10^{-6}$ & $\frac{373z^3 - 1139z^2 - 35117z + 10395}{37z^3 - 161z^2 - 1793z - 10395}$&\\
\cline{2-4}
& reduced & $\frac{504z^3 + 232z^2 - 616z + 222}{-168z^3 + 166z^2 + 397z - 114}$ & $\sm{7,-6,1,0}$   \\
\hline
$[0,1, 9, -15, -9, -3, -1, 15]$& $3.94828 \cdot 10^{-6}$ & $\frac{1701z^3 + 22095z^2 + 45819z + 5265}{-187z^3 - 1777z^2 + 5019z + 5265}$ & \\
\cline{2-4}
& reduced & $\frac{84z^3 + 664z^2 + 852z - 280}{-12z^3 - 5z^2 + 297z - 280}$ & $\sm{-3,-15,1,1}$   \\
\hline
\end{tabular}

\begin{tabular}{|l|c|c|c|}
\hline
Orbit & result & Interpolated Map  & conjugation\\
\hline
\multicolumn{4}{|c|}{Many Preperiodic Points}\\
\hline
$[0,1, 2, 3, 5, 4, -7, 6]$& 13 (6,2)  & $\frac{27z^3 - 298z^2 + 641z + 630}{9z^3 - 56z^2 - 83z + 630}$ &\\
\cline{2-4}
& reduced form& $\frac{130z^2 - 55z - 60}{-36z^3 - 50z^2 + 176z - 60}$ & $\sm{2,3,0,1}$  \\
\hline
$[0,-1, 1, -9, -15, 3, 9, 6]$& $13$ (6,2) &  $\frac{45z^3 + 243z^2 - 1998z - 810}{11z^3 + 23z^2 - 564z + 810}$ & \\
\cline{2-4}
& reduced form& $\frac{12z^3 + 94z^2 + 118z - 84}{33z^3 + 122z^2 - 43z - 42}$ & $\sm{3,3,0,1}$  \\
\hline
$[0,-1, -2, -4, -3, -10, -8, -5]$ & $13$ (2,6) & $\frac{8z^3 + 52z^2 + 34z - 220}{z^3 + 26z^2 + 140z + 220}$  &\\
\cline{2-4}
& reduced form& $\frac{24z^3 + 12^z2 - 39z + 3}{4z^3 + 28z^2 - 20z + 6}$ & $\sm{2,-4,0,1}$  \\
\hline
$[0,-1, -2, 5, 1, 4, -5, 3]$ & 12  (4,5) & $\frac{210z^3 - 775z^2 - 730z - 145}{83z^3 - 225z^2 - 363z + 145}$ &\\
\cline{2-4}
& reduced form& $\frac{180z^3 - 564z^2 - 48z + 332}{-270z^3 + 543z^2 - 169z - 254}$ & $\sm{1,-2,1,0}$  \\
\hline
$[0,1, 11, 5, -1, 2, -7, 3]$&12 (7,1) & $\frac{391z^3 - 2725z^2 + 4369z - 1155}{107z^3 - 965z^2 + 2093z - 1155}$ &\\
\cline{2-4}
& reduced form& $\frac{6z^3 + 7z^2 - 63z + 20}{9z^3 - 62z^2 + 33z + 20}$ & $\sm{4,-5,2,-1}$  \\
\hline
$[0,-1, -2, -3, -6, 2, 4, 1]$&12 (6,3)& $\frac{2z^2 - 26z + 12}{-x^3 - x^2 + 8z - 12}$ &\\
\cline{2-4}
& reduced form& $\frac{6z^3 + 8z^2 + 2z - 4}{3z^3 + 13z^2 + 2z}$ & $\sm{0,-2,1,0}$  \\
\hline
$[0,-1, -2, 5, 1, 4, -5, 3]$&12 (4,5)&$\frac{210z^3 - 775z^2 - 730z - 145}{83z^3 - 225z^2 - 363z + 145}$ &\\
\cline{2-4}
& reduced form& $\frac{180z^3 - 564z^2 - 48z + 332}{-270z^3 + 543z^2 - 169z - 254}$ & $\sm{1,-2,1,0}$  \\
\hline
\end{tabular}

\begin{tabular}{|l|c|c|c|}
\hline
Orbit & result & Interpolated Map   & conjugation \\
\hline
\hline
\multicolumn{4}{|c|}{Long Periodic Cycles}\\
\hline
\hline
$[0,1, 5, -9, -5, -3, -1, 3]$ & (0, 8) & $\frac{34z^3 - 21z^2 - 238z - 15}{-4z^3 + 19z^2 - 48z - 15}$ & \\
\cline{2-4}
&reduced & $\frac{38z^3 + 37z^2 - 39z - 24}{-8z^3 + 7z^2 - 11z - 12}$ & $\sm{2,1,0,1}$ \\
\hline
$[0,1, -1, 7, -5, -9, 3, -3]$ & (0, 8) & $\frac{59z^3 + 607z^2 + 1101z - 567}{-5z^3 - 113z^2 - 515z - 567}$ &  \\
\cline{2-4}
 & reduced & $\frac{24z^3 - 53z^2 - 45z + 20}{12z^3 - 77z^2 + x + 30}$ & $\sm{-5,1,1,-1}$ \\
\hline
$[0,-1, -3, -9, 3, 5, 15, 9]$ & (1, 8) & $\frac{153z^3 + 129z^2 - 14409z + 7695}{-z^3 + 223z^2 - 111z - 7695}$  & \\
\cline{2-4}
& reduced & $\frac{78z^3 + 72z^2 - 182z - 28}{-3z^3 + 107z^2 + 100z - 84}$ & $\sm{6,3,0,1}$  \\
\hline
$[0,1, 2, -1, -4, 5, -2, -5]$ & (0, 8) & $\frac{21z^3 + 76z^2 - 147z - 10}{-10z^3 - 29z^2 + 19z - 10}$ & \\
\cline{2-4}
& reduced & $\frac{16z^3 + 47z^2 - 3z - 90}{16z^3 + 63z^2 + 14z - 33}$  & $\sm{-1,-3,1,0}$  \\
\hline
$[0,-1, -7, 5, -5, 1, -3, -4]$ & (1, 8)&  $\frac{167z^3 + 1332z^2 + 2761z + 420}{-17z^3 - 276z^2 - 847z - 420}$  & \\
\cline{2-4}
& reduced & $\frac{28z^3 + 142z^2 - 138z + 108)}{-14z^3 - 147z^2 - 201z - 108}$ & $\sm{-1,-9,1,3}$  \\
\hline
$[0,1, -3, -7, -5, -9, 7, 3]$ & (1, 8) & $\frac{35z^3 + 175z^2 - 819z - 63}{3z^3 + 47z^2 + 237z - 63}$  & \\
\cline{2-4}
& reduced & $\frac{38z^3 + 27z^2 - 57z + 10}{12z^3 + 38z^2 + 38z - 16}$ & $\sm{4,-1,0,1}$  \\
\hline
\end{tabular}

\begin{tabular}{|l|c|c|c|}
\hline
Orbit & result & Interpolated Map & conjugation \\
\hline
\multicolumn{4}{|c|}{Long Preperiodic Tail}\\
\hline
$[0,-1, 2, 3, -3, 11, 1, 5]$&(8, 2) & $\frac{29z^3 + 538z^2 - 1837z - 330}{67z^3 + 14z^2 - 731z + 330}$  &\\
\cline{2-4}
& reduced form & $\frac{40z^3 - 61z^2 - 83z + 90}{-40z^3 + 15z^2 - 5z + 30}$ & $\sm{4,-5,2,-1}$ \\
\hline
$[0,1, 3, 5, 7, 21, 6, 9]$&(8, 2) & $\frac{144z^3 - 1491z^2 + 2322z + 945}{20z^3 - 139z^2 - 186z + 945}$ &\\
\cline{2-4}
&reduced form& $\frac{12z^3 - 40z^2 - 46z + 14}{-6z^3 + 27z^2 + 5z - 26}$ & $\sm{6,-3,2,1}$ \\
\hline
$[0,1, -1, -2, -7, -5, -3, 3]$&(7, 1) & $\frac{283z^3 + 683z^2 - 1131z - 315}{-9z^3 + 251z^2 + 553z - 315}$ &\\
\cline{2-4}
& reduced form& $\frac{137z^3 + 28z^2 - 203z + 38}{-9z^3 + 139z^2 + 6z - 76}$ &$\sm{2,-1,0,1}$ \\
\hline
$[0,-1, -2, -3, 7, 2, 3, -8]$&(7, 1) & $\frac{169z^3 - 423z^2 - 2188z + 2016}{4z^3 + 252z^2 + 38z - 2016}$ &\\
\cline{2-4}
&reduced form& $\frac{784z^3 + 18z^2 - 678z + 142}{748z^2 + 239z - 455}$ & $\sm{-8,1,1,-2}$  \\
\hline
$[0,1, 7, -1, -3, 5, -5, 3]$&(7, 1)& $\frac{14z^3 - 165z^2 + 214z + 105}{12z^3 - 45z^2 - 48z + 105}$ &\\
\cline{2-4}
&reduced form& $\frac{6z^3 + 51z^2 - 9z - 24}{-18z^3 + 7z^2 + 57z + 2}$ & $\sm{1,-2,1,0}$ \\
\hline
\end{tabular}

\end{center}
}

\begin{comment}
Code to analyze maps....
\begin{verbatim}
O =[-1, -2, -3, 7, 2, 3, -8]
F=orbit_to_map_rational([0] + O,3)
P = F.domain()
F.normalize_coordinates()
print(F.dehomogenize(1))

G,m=F.reduced_form()
G.dehomogenize(1), m

F.rational_preperiodic_graph()

hs = max([t.global_height() for t in F.sigma_invariants(1)])
hQ = F.canonical_height(P(0), error_bound=0.0000001)
hQ, hs, hQ/hs
\end{verbatim}

\subsubsection{Best Previously Known}
\begin{center}
\begin{tabular}{|l|c|c|}
\hline
Map & value & citation\\
\hline
\multicolumn{3}{|c|}{Long Periodic Cycles}\\
\hline
&&\\
\hline
\multicolumn{3}{|c|}{Long Preperiodic Tail}\\
\hline
&&\\\hline
\multicolumn{3}{|c|}{Many Preperiodic Points}\\
\hline
&&\\
\hline
\multicolumn{3}{|c|}{Small Height ratio}\\
\hline
&&\\
\hline
\end{tabular}
\end{center}
\end{comment}
%%%%%%%%%%%%%%%%%%%%%%%%%%%%%%%%%%%%%
\subsection{Degree 4: Polynomial}
{\small
\begin{center}
\begin{tabular}{|l|c|c|c|}
\hline
Orbit & result & Interpolated Map & conjugation \\
\hline
\hline
\multicolumn{4}{|c|}{Small Height ratio}\\
\hline
$[0,1, 4, 5, -1, 3]$ & $2.9015 \cdot 10^{-6}$ & $3/40z^4 - 14/15z^3 + 97/40z^2 + 43/30z + 1$ &\\
\cline{2-4}
& reduced form & $3/40z^4 - 14/15z^3 + 97/40z^2 + 43/30z + 1$ & $\sm{1,0,0,1}$\\
\hline
$[0,1, 7, 5, 8, 3]$ & $2.9568 \cdot 10^{-6}$ & $-11/840z^4 + 41/140z^3 - 2101/840z^2 + 1151/140z + 1$ & \\
\cline{2-4}
& reduced form & $-11/840z^4 - 3/140z^3 - 7/120z^2 - 207/140z + 4/7$ & $\sm{1,6,0,1}$ \\
\hline
$[0,1, -3, -2, 2, -4]$ & $3.0905\cdot 10^{-6}$ & $7/60x^4 + 11/20x^3 - 29/30x^2 - 37/10x + 1$ &\\
\cline{2-4}
& reduced form & $7/60x^4 - 79/60x^3 + 109/30x^2 + 17/30x - 1$ & $[1,-4,0,1]$ \\
\hline
$[0,1, 3, -3, -4, 8]$& $3.6227\cdot 10^{-6}$ & $3/56x^4 - 29/84x^3 - 55/56x^2 + 275/84x + 1$ &\\
\cline{2-4}
& reduced form & $3/56x^4 - 11/84x^3 - 95/56x^2 + 41/84x + 2$ & $[1,1,0,1]$  \\
\hline
$[0,1, -2, 8, 7, 5]$& $3.8889\dot 10^{-6}$ & $-3/280x^4 + 47/420x^3 + 87/280x^2 - 1433/420x + 1$ &\\
\cline{2-4}
& reduced form  & $-3/280x^4 - 97/420x^3 - 313/280x^2 + 463/420x - 1$ & $[1,8,0,1]$  \\
%\hline
%$[1, 5, -2, -7, -9]$& $3.8646\cdot 10^{-6}$ & $-1/120x^4 - 11/84x^3 - 89/840x^2 + 1783/420x + 1$&2000\\
%& reduced & $-1/120x^4 + 43/420x^3 + 163/840x^2 - 877/420x - 2$ & $[1,-7,0,1]$ & \\
%\hline
%$[1, -7, -1, -2, -9]$& $3.9572\cdot 10^{-6}$ & $-71/840x^4 - 187/140x^3 - 4549/840x^2 - 163/140x + 1$ &2000\\
%& reduced & same & id & \\
\hline
\hline
\multicolumn{4}{|c|}{Many Preperiodic Points}\\
\hline
$[0,1, 10, 4, 8, 3]$ & 11 $(1, 6)$ &$-1/40z^4 + 37/60z^3 - 199/40z^2 + 803/60z + 1$ &\\
\cline{2-4}
& reduced form & $-1/40z^4 + 7/60z^3 + 21/40z^2 - 157/60z$ & $\sm{1,5,0,1}$ \\
\hline
$[0,1, 6, -3, -6, -1]$ & 10 $(3, 3)$ & $-1/60z^4 - 2/15z^3 + 31/60z^2 + 139/30z + 1$ &\\
\cline{2-4}
& reduced form & $-1/60z^4 + 2/15z^3 + 31/60z^2 - 49/30z - 1$ & $\sm{1,-4,0,1}$  \\
\hline
\hline
\multicolumn{4}{|c|}{Long Periodic Cycles}\\
\hline
\hline
$[0,1, 5, 4, 6, 2]$ & (1, 6) & $1/120z^4 - 1/60z^3 - 121/120z^2 + 301/60z + 1$ &\\
\cline{2-4}
 & reduced form & $1/120z^4 + 7/60z^3 - 49/120z^2 - 103/60z + 2$ & $\sm{1,4,0,1}$ \\
\hline
$[0,1, 3, -1, -2, -5]$ & (1, 6) & $-1/120z^4 - 11/60z^3 - 59/120z^2 + 161/60z + 1$ &\\
\cline{2-4}
& reduced form & $-1/120z^4 - 1/4z^3 - 43/24z^2 - 7/4z + 4/5$ & $\sm{1,2,0,1}$ \\
\hline
$[0,1, 10, 4, 8, 3]$ & (1, 6) & $-1/40z^4 + 37/60z^3 - 199/40z^2 + 803/60z + 1$ &\\
\cline{2-4}
&reduced form & $-1/40z^4 + 7/60z^3 + 21/40z^2 - 157/60z$ & $\sm{1,5,0,1}$ \\
\hline
$[0,1, 4, -1, -2, 3]$ & (0, 6) & $11/60z^4 - 29/30z^3 - 11/60z^2 + 119/30z + 1$ &\\
\cline{2-4}
&reduced form & $11/60z^4 + 1/2z^3 - 19/12z^2 - 5/2z + 7/5$ & $\sm{1,2,0,1}$\\
\hline
\hline
\multicolumn{4}{|c|}{Long Preperiodic Tail}\\
\hline
$[0,1, -1, 2, -3, 4]$ & $(5, 1)$ 7 & $7/120z^4 + 1/20z^3 - 67/120z^2 - 31/20z + 1$ &\\
\cline{2-4}
& reduced form  & $7/120z^4 + 1/20z^3 - 67/120z^2 - 31/20z + 1$ & $\sm{1,0,0,1}$\\
\hline
$[0,1, -1, 2, -2, 3]$ & $(5, 1)$ 7& $1/8z^4 + 1/12z^3 - 5/8z^2 - 19/12z + 1$&\\
\cline{2-4}
& reduced form & $1/8z^4 + 1/12z^3 - 5/8z^2 - 19/12z + 1$ & $\sm{1,0,0,1}$\\
\hline
$[0,1,-3,-2,2,-5]$ & $(5,1)$ 8 & $11/120z^4 + 9/20z^3 - 119/120z^2 - 71/20z + 1$ &\\
\cline{2-4}
& reduced form & $11/120z^4 - 61/60z^3 + 289/120z^2 + 151/60z - 2$ &$\sm{1,-4,0,1}$ \\
\hline
$[0,1, 5, -1, 2, 9]$ & $(5, 1)$ 7 & $7/90z^4 - 89/90z^3 + 109/45z^2 + 112/45z + 1$ &\\
\cline{2-4}
& reduced form & $7/90z^4 - 89/90z^3 + 109/45z^2 + 112/45z + 1$ & $\sm{1,0,0,1}$\\
\hline
\end{tabular}
\end{center}
}

\begin{comment}
Code to analyze maps....
\begin{verbatim}
# what to do about normal form???
P.<x,y>=ProjectiveSpace(QQ,1)
F=orbit_to_map_polynomial([1, 5, 4, 6, 2],4)
print(F.dehomogenize(1))
f=DynamicalSystem(list(F), domain=P)
g=f.normal_form()
print(g.dehomogenize(1))
f.reduced_form()

gr=f.rational_preperiodic_graph()
gr.show()
print(gr.num_verts())

print(f.rational_preimages(P(0)))

hF = max([t.global_height() for t in f[0].coefficients()])
hs = max([t.global_height() for t in f.sigma_invariants(1)])
hQ=f.canonical_height(P(0), error_bound=0.0000001)
hQ, hF, hs, hQ/hF, hQ/hs
\end{verbatim}

\subsubsection{Best Previously Known}
\begin{center}
\begin{tabular}{|l|c|c|}
\hline
Map & value & citation\\
\hline
\multicolumn{3}{|c|}{Long Periodic Cycles}\\
\hline
&&\\
\hline
\multicolumn{3}{|c|}{Long Preperiodic Tail}\\
\hline
&&\\\hline
\multicolumn{3}{|c|}{Many Preperiodic Points}\\
\hline
$1/24x^4 - 11/12x^3 + 167/24x^2 - 253/12x + 23$& 11 (2,2) &Doyle-Hyde\\
\hline
\multicolumn{3}{|c|}{Small Height ratio}\\
\hline
&&\\
\hline
\end{tabular}
\end{center}
\end{comment}

\subsection{Degree 4: Rational Functions}
{\small
\begin{center}
\begin{tabular}{|l|c|c|c|}
\hline
Orbit & result & Interpolated Map  & conjugation \\
\hline
\hline
\multicolumn{4}{|c|}{Small Height ratio}\\
\hline
$[0,-1, -4, -3, -9, 2, -6, -2, -5, 1]$ & $2.1843\cdot 10^{-8}$ & $\frac{49x^4 + 844x^3 + 4499x^2 + 7364x - 660}{-10x^4 - 174x^3 - 866x^2 - 1122x + 660}$ &\\
\cline{2-4}
&reduced form & $\frac{9x^4 + 4x^3 + 123x^2 - 604x - 132}{-10x^4 - 14x^3 + 262x^2 + 14x - 132}$ & $\sm{1,-4,0,1}$ \\
\hline
$[0,1, 8, 5, -4, -1, -3, -2, 2, 4]$ & $2.6346\cdot 10^{-8}$ & $\frac{(91x^4 + 428x^3 - 2133x^2 - 4922x - 184}{60x^4 - 236x^3 - 674x^2 + 194x - 184}$ &\\
\cline{2-4}
&reduced form & $\frac{211x^4 - 1732x^3 + 1847x^2 + 2110x - 1680}{60x^4 - 716x^3 + 2182x^2 - 1862x - 420}$ & $\sm{1,-2,0,1}$ \\
\hline
$[0,-1, -4, -5, 3, -9, -6, -3, -2, 1]$ & $2.9735 \cdot 10^{-8}$ & $\frac{182x^4 + 1323x^3 + 2023x^2 + 717x + 1755}{-15x^4 - 36x^3 - 30x^2 - 1284x - 1755}$ & \\
\cline{2-4}
&reduced form & $\frac{334x^4 + 619x^3 - 433x^2 - 340x + 180}{-60x^4 + 48x^3 - 12x^2 - 636x - 120}$ & $\sm{2,-1,0,1}$ \\
\hline
$[0,-1, 1, -2, 5, -4, 3, -5, -9, -3]$ & $3.2505\cdot 10^{-8}$ & $\frac{52x^4 + 503x^3 + 347x^2 - 4187x - 2475}{-11x^4 - 112x^3 - 220x^2 + 748x + 2475}$ &  \\
\cline{2-4}
&reduced form & $\frac{38x^4 - 61x^3 - 395x^2 + 598x + 624}{-44x^4 + 40x^3 + 194x^2 + 116x + 96}$ & $\sm{2,-3,0,1}$ \\
\hline
$[0,-1, -4, -5, -2, -3, 5, 3, 2, 1]$ & $3.4953 \cdot 10^{-8}$ & $\frac{27x^4 + 220x^3 + 166x^2 - 628x - 25}{8x^4 + 43x^3 + 135x^2 + 269x + 25}$ & \\
\cline{2-4}
&reduced form & $\frac{70x^4 + 123x^3 - 139x^2 - 78x + 54}{32x^4 + 22x^3 + 54x^2 + 48x - 36}$ & $\sm{2,-1,0,1}$ \\
%$ [[-1, -9, -7, 6, 5, 2, -3, 3, -2], 3.47796040632976e-8, 0],$&&&&\\
%$ [[1, -5, 3, -2, -3, 6, 5, 4, -1], 3.63671009696441e-8, 0],$&&&&\\
%$ [[-1, -5, -6, -7, -8, -12, -4, -2, -9], 3.66574866822840e-8, 0],$&&&&\\
%$ [[-1, 9, -7, -3, -9, 1, 6, -6, -4], 4.48712067618227e-8, 0],$&&&&\\
%$ [[-1, 5, -3, 11, 4, -5, 7, 6, -4], 4.62780992132954e-8, 0],$&&&&\\
%$[[1, -1, 5, -4, -11, -3, -2, 3, -5], 4.80293376753452e-8, 0]$&&&&\\
%$[[-1, 9, -2, 7, 3, -3, 5, -7, 1], 4.85583864697068e-8, 0]$&&&&\\
%$[[1, -1, 5, -5, -7, -4, -9, -2, -3], 4.96171482909717e-8, 0]$&&&&\\
%$[1, 8, -6, 2, -2, 4, 3, -16, -3];6.71308197247157e-8;0;99438;2;$&&&&\\
%$[1, -7, -6, -2, -1, -8, 2, -11, -4];6.73376328384607e-8;0;99438;3;$&&&&\\
%$[-1, -6, 15, -15, 1, -3, 3, -7, -9];7.46685662292143e-8;0;99438;3;$&&&&\\
%$[1, 6, -7, 4, -1, 10, -8, -3, 11];7.49039890105353e-8;0;99438;3;$&&&&\\
%$[-1, -6, -3, 10, 1, -10, -4, 15, -2];7.81093452591267e-8;0;99438;0;$&&&&\\
%$[1, 2, 5, 9, 15, 17, 21, -19, -1];7.89946146003784e-8;0;99438;7;$&&&&\\
%$[-1, -9, -3, -11, -12, 4, -4, 9, 1];7.90346079785275e-8;0;99438;2;$&&&&\\
%$[-1, -7, 23, -13, 14, -23, 15, -21, -5];8.07064011842145e-8;0;99438;2;$&&&&\\
\hline
\hline
\multicolumn{4}{|c|}{Many Preperiodic Points}\\
\hline
$[0,1, 3, 2, 5, -2, -3, -4, -5, -9]$& 11 $(7, 3)$ & $\frac{61x^4 + 73x^3 - 1441x^2 + 1007x + 540}{-13x^4 - 15x^3 + 161x^2 - 593x + 540}$ & \\
\cline{2-4}
&reduced form & $\frac{148x^4 + 384x^3 - 447x^2 - 261x + 20}{-52x^4 - 134x^3 + 38x^2 - 184x + 20}$ & $\sm{2,1,0,1}$ \\
\hline
%$[-1, -3, -5, -7, 1, 3, 2, 7, 5];-10;(1, 9);99053;1;$&&&&\\
%$[1, -9, 7, -15, 5, -7, -1, 9, -5];-10;(8, 2);99053;44;$&&&&\\
%$[1, -7, -1, -2, -5, 5, 4, 3, 2];-10;(8, 2);99064;49;$&&&&\\
$[1, -1, 5, 7, 3, -5, -4, -3, -2]$ & 11 (1, 10) & $\frac{39x^4 - 67x^3 + 991x^2 + 2107x - 910}{13x^4 + 115x^3 - 375x^2 - 1003x - 910}$ & \\
\cline{2-4}
&reduced form & $\frac{12x^4 + 5x^3 + 321x^2 + 126x - 104}{-18x^4 + 6x^3 + 98x^2 + 80x + 104}$ & $\sm{-1,-4,1,0}$\\
\hline
$[0,1, -2, 7, -1, 3, -3, -7, 2, -5]$ & $11, (1, 9)$ & $\frac{13x^4 - 80x^3 - 136x^2 + 1580x - 609}{2x^4 - 45x^3 + 43x^2 + 225x - 609}$ & \\
&reduced form & $\frac{22x^4 + 9x^3 - 109x^2 + 234x + 144}{8x^4 - 74x^3 - 80x^2 + 92x - 96}$ & $\sm{2,1,0,1}$ \\
\hline
\hline
\multicolumn{4}{|c|}{Long Periodic Cycles}\\
\hline
\hline
$ [0,-1, -2, -4, 3, 6, 5, 4, -5, 1]$ & $(0, 10)$ & $\frac{18x^4 - 143x^3 - 828x^2 + 5963x - 5010}{-14x^4 + 158x^3 - 391x^2 - 1373x + 5010}$ &\\
\cline{2-4}
&reduced form & $\frac{88x^4 + 827x^3 + 332x^2 - 1877x + 630}{-14x^4 - 122x^3 - 121x^2 - 433x - 630}$ & $\sm{1,5,0,1}$  \\
\hline
%$ [[-1, -3, -5, -7, 1, 3, 2, 7, 5], -14, (1, 9)],$&&&&\\
$ [0,-1, -7, -9, 7, -5, -3, 5, 3, 1]$ & $(0, 10)$ & $\frac{135*x^4 - 948*x^3 - 12546*x^2 - 27276*x + 40635}{161*x^4 + 1068*x^3 - 5798*x^2 - 39276*x - 40635}$ &\\
\cline{2-4}
&reduced form & $\frac{309*x^4 - 1290*x^3 - 2109*x^2 + 1794*x + 2304}{161*x^4 - 432*x^3 - 1679*x^2 + 870*x + 576}$ & $\sm{2,-3,0,1}$ \\
\hline
$[0,1, 3, 5, 7, -1, -3, -2, -7, -5]$ & $(1, 9)$ & $\frac{431x^4 + 1042x^3 - 19244x^2 - 67138x - 19635}{59x^4 + 382x^3 - 2456x^2 - 13198x - 19635}$&\\
\cline{2-4}
&reduced form & $\frac{245x^4 - 134x^3 - 2879x^2 - 2164x + 576}{59x^4 + 73x^3 - 812x^2 - 922x - 576}$ & $\sm{2,-1,0,1}$ \\
\hline
$[0,1, -3, -5, -1, 3, 2, 5, 4, 9]$ & $(1, 9)$ & $\frac{38x^4 - 65x^3 - 2759x^2 + 9929x - 3975}{52x^4 - 479x^3 + 635x^2 + 2711x - 3975}$ &\\
\cline{2-4}
&reduced form & $\frac{118x^4 + 22x^3 - 328x^2 + 236x + 72}{-104x^4 - 145x^3 + 434x^2 + 199x - 144}$ & $\sm{2,3,0,1}$  \\
\hline
%$ [[-1, -3, -5, -7, 1, 3, 2, 7, 5], -14, (1, 9)],$&&&&\\
%$ [[1, 7, 5, -7, -9, -3, 2, 3, -1], -14, (1, 9)],$&&&&\\
\hline
\hline
\multicolumn{4}{|c|}{Long Preperiodic Tail}\\
\hline
$[0,-1, 1, 9, 5, -11, -9, -7, -5, -3]$ & $(9, 1)$ & $\frac{223x^4 - 294x^3 - 5988x^2 - 10066x - 1155}{x^4 + 254x^3 - 396x^2 - 2934x + 1155}$ &\\
\cline{2-4}
&reduced form & $\frac{222x^4 + 85x^3 - 369x^2 - 298x - 60}{4x^4 + 258x^3 + 93x^2 - 185x - 30}$ & $\sm{4,1,0,1}$  \\
\hline
$[0,-1, -4, 4, 2, -2, 1, -3, 3, 6]$ & $(8, 2)$ & $\frac{40x^4 - 153x^3 - 67x^2 + 558x + 72}{8x^4 - 20x^3 + 34x^2 - 100x - 72}$ &\\
\cline{2-4}
&reduced form & $\frac{40x^4 - 153x^3 - 67x^2 + 558x + 72}{8x^4 - 20x^3 + 34x^2 - 100x - 72}$ & $\sm{1,0,0,1}$\\
\hline
$[0,-1, 3, -3, 4, -5, 1, -2, 2, -6]$ & $(8, 2)$ & $\frac{21x^4 + 72x^3 - 153x^2 - 264x - 180}{2x^4 - 22x^3 - 76x^2 + 168x + 180}$ &\\
\cline{2-4}
&reduced form & $\frac{21x^4 + 72x^3 - 153x^2 - 264x - 180}{2x^4 - 22x^3 - 76x^2 + 168x + 180}$ & $\sm{1,0,0,1}$\\
\hline
%$[1, 2, -2, -1, -3, 3, -5, 7, -7];-18;(8, 2);99428;40;$&&&&\\
%$[1, 3, -5, 7, -1, -4, 2, 5, -7];-18;(8, 2);99429;11;$&&&&\\
%$ [[-1, -3, 2, 1, 3, 5, 7, 6, -7], -18, (8, 2)],$&&&&\\
%$ [[-1, 1, -2, -3, 5, 7, 3, -7, -5], -18, (8, 2)],$&&&&\\
$[0,1, -5, -7, -1, 3, 9, 7, -3, 5]$ & $(9, 1)$ & $\frac{91x^4 + 418x^3 - 4696x^2 - 16738x + 5565}{27x^4 - 30x^3 - 1176x^2 - 1314x + 5565}$ & \\
\cline{2-4}
&reduced form & $\frac{32x^4 + 88x^3 - 56x^2 - 163x - 36}{54x^4 + 39x^3 - 138x^2 - 114x + 24}$ & $\sm{4,1,0,1}$  \\
\hline
\end{tabular}
\end{center}
}
\begin{comment}
\begin{verbatim}
params['map_type'] = 'rational'
params['degree'] = 4
params['population'] = 1000
params['generations'] = 1000
params['survival'] = 0.15
params['reset_survival'] = 0.02
params['reset_interval'] = 50
params['normalize_orbit'] = True
params['bound'] = 20
params['mixing_method'] = 'permutation'
params['mutation_rate'] = 0.05
params['mutation_method'] = 'all'
params['target'] = 'preperiodic'
params['orbit_target'] = 11
params['orbit_weights'] = (5,1)
\end{verbatim}
\end{comment}

%%%%%%%%%%%%%%%%%%%%%%%%%%%%%%%%%%%%%
\subsection{Degree 5: Polynomial}
{\small
\begin{center}
\begin{tabular}{|l|c|c|c|}
\hline
Orbit & result & Interpolated Map & conjugation \\
\hline
\hline
\multicolumn{4}{|c|}{Small Height ratio}\\
\hline
$[0,-1, -2, -6, -4, -5, -3]$ & $9.1519\cdot 10^{-9}$ & $-1/40x^5 - 1/2x^4 - 27/8x^3 - 17/2x^2 - 23/5x - 1$& \\
reduced form && $-2/5x^5 + 5/2x^3 - 13/5x - 1/2$ & $\sm{2,-4,0,1}$ \\
\hline
%;1.24727487024804e-8;0$&&&&\\
$[0,-1, -4, -3, -6, -5, -2]$ & $1.6874\cdot 10^{-8}$ & $-1/30x^5 - 7/12x^4 - 10/3x^3 - 77/12x^2 - 19/30x - 1$ &\\
\cline{2-4}
 &reduced form& $-1/30x^5 + 1/12x^4 + 2/3x^3 - 13/12x^2 - 79/30x + 1$ & $\sm{1,-4,0,1}$ \\
\hline
%$[1, 4, 3, 6, 5, 2];1.70586107374925e-8;0;101927;2;$&&&&\\
$[0,-1, -4, -11, -8, -15, -12]$ & $2.8227\cdot 10^{-8}$ & \scalebox{0.7}{$13/18480x^5 + 71/3696x^4 + 587/3696x^3+ 2281/3696x^2 + 16063/4620x - 1$} & \\
\cline{2-4}
 &reduced form& \scalebox{0.7}{$13/1155x^5 - 23/231x^4 + 137/924x^3+ 248/231x^2 - 9137/4620x - 353/154$} & $\sm{2,-9,0,1}$ \\
\hline
$[0,-1, -8, -9, -10, -5, -3]$ & $3.1906 \cdot 10^{-8}$ & \scalebox{0.9}{$-1/840x^5 + 1/84x^4 + 101/168x^3 + 395/84x^2 + 389/35x - 1$} &\\
\cline{2-4}
 &reduced form& \scalebox{0.9}{$-2/105x^5 + 8/21x^4 - 19/42x^3 - 41/21x^2 + 113/70x + 17/1$} & $\sm{2,-6,0,1}$ \\
%$[-1, -5, -8, -7, -9, 1];4.39342795397197e-8;0$&&&&\\
%$[-1, -5, -3, -4, -2, -6];5.12946858565901e-8;0$&&&&\\
%$[-1, -7, -6, -3, -8, 1];5.32134184790335e-8;0$&&&&\\
%$[1, -3, -8, -1, -4, -6];6.24297764984559e-8;0$&&&&\\
\hline
\end{tabular}

\begin{tabular}{|l|c|c|c|}
\hline
Orbit & result & Interpolated Map & conjugation \\
\hline
\hline
\multicolumn{4}{|c|}{Many Preperiodic Points}\\
\hline
$[0,1, 3, 5, -7, -9, -11]$ & $12 (1, 6)$& $-1/840x^5 - 1/56x^4 - 3/56x^3 + 9/56x^2 + 803/420x + 1$ &\\
\cline{2-4}
 &reduced form& $-2/105x^5 + 3/14x^3 + 199/210x$ & $\sm{2,-3,0,1}$ \\
\hline
$[0,-1, 2, -5, 1, -4, 3]$ & $11 (1, 6)$& $1/840x^5 + 1/28x^4 + 43/168x^3 - 1/28x^2 - 114/35x - 1$&\\
\cline{2-4}
 &reduced form& $2/105x^5 - 29/42x^3 + x^2 + 117/70x + 5/2$ & $\sm{2,-6,0,1}$ \\
\hline
$[0,-1, 3, -7, -6, -4, 5]$ & $11 (5, 2)$ & $-1/315x^5 - 1/126x^4 + 2/9x^3 + 11/63x^2 - 851/210x - 1$&\\
\cline{2-4}
&reduced form& $-1/315x^5 + 1/126x^4 + 2/9x^3 - 32/63x^2 - 781/210x + 4$ & $\sm{1,-1,0,1}$ \\
\hline
$[0,1, -2, 2, -3, -1, 3]$ & $11 (3, 7)$ & $-1/60x^5 + 1/24x^4 + 1/2x^3 - 13/24x^2 - 179/60x + 1$ &\\
\cline{2-4}
 &reduced form& $-1/60x^5 + 7/24x^4 - 3/2x^3 + 41/24x^2 + 151/60x + 2$ & $\sm{1,-3,0,1}$ \\
\hline
%$[-1, 2, -2, 3, 1, -3];-11;(3, 7)$&&&&\\
$[0,-1, -2, -4, -3, -5, 3]$ & $11 (4, 3)$ & $1/120x^5 + 1/12x^4 - 1/24x^3 - 13/12x^2 + 1/30x - 1$&\\
reduced form && $2/15x^5 - 3/2x^3 + x^2 + 28/15x - 1$ & $\sm{2,-2,0,1}$ \\
%$[1, -3, -5, -2, -4, -1];-10;(0, 7)$&&&&\\
%$[-1, -8, -3, -7, -6, 1];-10;(0, 8)$&&&&\\
\hline
\end{tabular}

\begin{tabular}{|l|c|c|c|}
\hline
Orbit & result & Interpolated Map & conjugation \\
\hline
\hline
\multicolumn{4}{|c|}{Long Periodic Cycles}\\
\hline
\hline
$[0,1, 7, 5, -1, 3, 8]$ & $(0, 8)$ & $-43/1680x^5 + 157/336x^4 - 125/48x^3 + 1187/336x^2 + 3889/840x + 1$ &\\
\cline{2-4}
 &reduced form& $-43/105x^5 + 2/3x^4 + 235/84x^3 - 19/6x^2 - 1843/420x + 5/2$ & $\sm{2,3,0,1}$ \\
\hline
$[0,1, 8, 3, 7, 6, -1]$ & $(0, 8)$ & $-47/1260x^5 + 47/72x^4 - 218/63x^3 + 313/72x^2 + 2309/420x + 1$ &\\
\cline{2-4}
 &reduced form& $-47/1260x^5 + 47/72x^4 - 218/63x^3 + 313/72x^2 + 2309/420x + 1$ & $\sm{1,0,0,1}$\\
\hline
$[0,1, -1, 9, 7, 3, 6]$ & $(0, 8)$& $11/1680x^5 - 25/336x^4 - 15/112x^3 + 1033/336x^2 - 4093/840x + 1$ &\\
\cline{2-4}
 &reduced form& $11/105x^5 + 5/7x^4 + 5/84x^3 - 107/28x^2 - 32/35x + 20/7$ & $\sm{2,5,0,1}$ \\
\hline
\end{tabular}
\begin{tabular}{|l|c|c|c|}
\hline
Orbit & result & Interpolated Map & conjugation \\
\hline
\hline
\multicolumn{4}{|c|}{Long Preperiodic Tail}\\
\hline
$[0,-1, -2, 3, 5, 6, 7]$ & $(6, 2)$ & $-13/2520x^5 + 23/252x^4 - 37/72x^3 + 187/252x^2 + 247/105x - 1$&\\
\cline{2-4}
 &reduced form& $-26/315x^5 - 2/21x^4 + 61/126x^3 + 1/14x^2 + 121/315x + 31/42$ & $\sm{2,4,0,1}$ \\
\hline
$[0,-1, 6, 4, 8, 1, 7]$ & $(6, 1)$ & \scalebox{0.9}{$-181/5040x^5 + 167/252x^4 - 565/144x^3 + 1723/252x^2 + 1873/420x - 1$} &\\
\cline{2-4}
 &reduced form& failed &\\
\hline
$[0,1, 3, 2, 4, -1, 6]$ & $(6, 1)$ & $-1/60x^5 + 1/3x^4 - 23/12x^3 + 19/6x^2 + 13/30x + 1$&\\
\cline{2-4}
 &reduced form& $-1/60x^5 + 1/6x^4 + 1/12x^3 - 5/3x^2 - 17/30x + 2$ & $\sm{1,2,0,1}$ \\
\hline
\end{tabular}

\end{center}
}

\begin{comment}
\subsubsection{Best Previously Known}
\begin{center}
\begin{tabular}{|l|c|c|}
\hline
Map & value & citation\\
\hline
\multicolumn{3}{|c|}{Long Periodic Cycles}\\
\hline
&&\\
\hline
\multicolumn{3}{|c|}{Long Preperiodic Tail}\\
\hline
&&\\\hline
\multicolumn{3}{|c|}{Many Preperiodic Points}\\
\hline
$1/120x^5 - 7/24x^4 + 89/24x^3 - 497/24x^2 + 2897/60x - 30$& 14 (4,1) &Doyle-Hyde\\
\hline
\multicolumn{3}{|c|}{Small Height ratio}\\
\hline
&&\\
\hline
\end{tabular}
\end{center}
\end{comment}

\subsection{Degree 5: Rational Functions}
{\small
\begin{center}
\begin{tabular}{|l||}
\hline
Orbit/Map  \\
\hline
\multicolumn{1}{|c|}{Small Height ratio}\\
\hline
$[0,1, -1, 7, -17, 13, -27, -8, -2, 14, -29, 4]$ \quad $3.6941 \cdot 10^{-10}$ \\
$\frac{88758966490337x^5 + 2981527184876307x^4 + 11066877636602289x^3 - 75987409403569835x^2 - 1288812439794810x + 17419787236064304}{1157006057963x^5 - 50301712572351x^4 - 1971337131131061x^3 + 814941375750391x^2 + 29505024045162162x + 17419787236064304}$ \\
\hline
$0,[1, -14, 22, 14, -20, 7, -11, -2, 16, 25, -1]$ \quad $3.7088 \cdot 10^{-10}$ \\
$\frac{360307694x^5 - 24850051494x^4 + 445661486158x^3 + 1287816130714x^2 - 68914533709032x + 4016025536600}{224318693x^5 - 5595423774x^4 - 29944225157x^3 + 974275081474x^2 - 441448123596x + 4016025536600}$ \\
\hline
$[0,-1, -91, -15, 1, 17, 6, -19, 23, -5, 2, -6]$ \quad $3.7574 \cdot 10^{-10}$ \\
\scalebox{0.7}{$\frac{410868324212067295x^5 - 214908806651361490x^4 - 123140623561248668490x^3 - 423945106692483245240x^2 + 966459488711406110795x - 3269779023770435670}{-31058519941648543x^5 - 315256851535218746x^4 + 8446975942617674502x^3 + 16273986690252479876x^2 - 3156194582136401159x + 3269779023770435670}$} \\
\hline
\hline
\end{tabular}
\end{center}
}

%%%%%%%%%%%%%%%%%%%%%%%%%%%%%%%%%%%%%
\subsection{Degree 6: Polynomial}
{\small
\begin{center}
\begin{tabular}{|l|c|c|c|}
\hline
Orbit & result & map & conjugation \\
\hline
\hline
\multicolumn{4}{|c|}{Small Height ratio}\\
\hline
$ [0,1, -1, -4, -7, -3, -6, -9]$ & $1.8372 \cdot 10^{-9}$ & $-29/5040x^6 - 71/560x^5 - 985/1008x^4$ &\\
&&$- 341/112x^3 - 6343/2520x^2 + 327/70x + 1$&\\
\cline{2-4}
 &reduced form& $-29/5040x^6 - 13/560x^5 + 149/1008x^4$ & $\sm{1,-3,0,1}$ \\
&&$+ 127/336x^3 - 809/1260x^2 + 61/420x - 3$&\\
\hline
$[0,-1, -8, -9, -7, -10, -2, -4]$ & $2.0560 \cdot 10^{-9}$ & $1/420x^6 + 137/1680x^5 + 181/168x^4$ &\\
&&$+ 753/112x^3 + 16193/840x^2 + 2158/105x - 1$&\\
\hline
\end{tabular}
\end{center}
}
{\small
\begin{center}
\begin{tabular}{|l|c|c|c|}
\hline
Orbit & result & map & conjugation\\
\hline
\hline
\multicolumn{4}{|c|}{Many Preperiodic Points}\\
\hline
$[0,1, -9, -2, -8, -3, -7, -1]$ &12 $(5, 4)$& $-11/5040x^6 - 39/560x^5 - 811/1008x^4$&\\
&&$- 1327/336x^3 - 16867/2520x^2 + 319/210x + 1$&\\
\cline{2-4}
 &reduced form& $-11/5040x^6 - 17/560x^5 - 55/1008x^4$ & $\sm{1,-3,0,1}$ \\
 &&$+ 69/112x^3 + 3923/2520x^2 - 108/35x - 4$&\\
\hline
$[0,1, -3, -4, -5, -1, 2, -7]$ & 12 $(5, 3)$ & $-1/504x^6 - 1/140x^5 + 59/504x^4$  &\\
&&$+ 9/28x^3 - 407/252x^2 - 197/70x + 1$&\\
\hline
\hline
\end{tabular}
\end{center}
}

{\small
\begin{center}
\begin{tabular}{|l|c|c|c|}
\hline
Orbit & result & map & conjugation\\
\hline
\hline
\multicolumn{4}{|c|}{Long Periodic Cycles}\\
\hline
\hline
$[0,1, -8, -5, -1, 2, -7, -4]$ & $(1, 8)$ & $37/5040x^6 + 233/1680x^5 + 113/144x^4$ &\\
&& $+ 223/336x^3 - 3019/630x^2 - 2437/420x + 1$&\\
\cline{2-4}
 &reduced&$37/5040x^6 + 53/560x^5 + 29/144x^4$ & $\sm{1,-1,0,1}$ \\
&form& $- 415/336x^3 - 4223/1260x^2 + 197/60x + 3$&\\
\hline
$[0,-1, 4, 1, -6, -3, -5, -2]$ &$(1, 8)$&$-19/1680x^6 - 47/560x^5 + 61/336x^4$ &\\
&&$+ 617/336x^3 - 143/840x^2 - 709/105x - 1$&\\
\hline
\hline
\end{tabular}
\end{center}
}
{\small
\begin{center}
\begin{tabular}{|l|c|c|c|}
\hline
Orbit & result & map & conjugation\\
\hline
\hline
\multicolumn{4}{|c|}{Long Preperiodic Tail}\\
\hline
$[0,1, 2, 7, 6, 5, 8, 9]$ & $(7, 2)$ & $-1/280x^6 + 59/840x^5 - 67/168x^4$&\\
&&$+ 25/168x^3 + 1429/420x^2 - 233/105x + 1$&\\
\cline{2-4}
 &reduced form& $-1/280x^6 - 7/120x^5 - 37/168x^4$ & $\sm{1,6,0,1}$ \\
&&$+ 73/168x^3 + 467/210x^2 - 289/210x - 1$&\\
\hline
$[0,1, 3, 7, 4, 2, 8, 9]$ & $(7, 2)$ & $13/3360x^6 - 59/480x^5 + 965/672x^4$ &\\
&&$- 4967/672x^3 + 3683/240x^2 - 509/70x + 1$&\\
\hline
\hline
\end{tabular}
\end{center}
}

\begin{comment}
\subsubsection{Best Previously Known}
\begin{center}
\begin{tabular}{|l|c|c|}
\hline
Map & value & citation\\
\hline
\multicolumn{3}{|c|}{Long Periodic Cycles}\\
\hline
&&\\
\hline
\multicolumn{3}{|c|}{Long Preperiodic Tail}\\
\hline
&&\\\hline
\multicolumn{3}{|c|}{Many Preperiodic Points}\\
\hline
$1/720x^6 - 1/16x^5 + 155/144x^4 - 425/48x^3 + 1594/45x^2 - 751/12x + 43$& 15 (3,2) &Doyle-Hyde\\
\hline
\multicolumn{3}{|c|}{Small Height ratio}\\
\hline
&&\\
\hline
\end{tabular}
\end{center}
\end{comment}

%%%%%%%%%%%%%%%%%%%%%%%%%%%%%%%%%%%%%
\subsection{Degree 7: Polynomial}
{\small
\begin{center}
\begin{tabular}{|l|c|c|c|}
\hline
Orbit & result & Interpolated Map & conjugation \\
\hline
\hline
\multicolumn{4}{|c|}{Small Height ratio}\\
\hline
$[-1, 1, 11, 2, 10, $ & $1.0564 \cdot 10^{-10}$ & $-43/33264x^7 + 691/15120x^6 - 100973/166320x^5 $ &\\
$8, 9, 3]$&&$+ 11099/3024x^4- 38399/4158x^3 + 24827/7560x^2$&\\
&&$ + 205729/13860x - 1$&\\
 &reduced form& unknown &\\
\hline
$[1, 5, -1, 14, 2, $ & $1.4693 \cdot 10^{-10}$ & $-71/8648640x^7 + 11383/8648640x^6 - 390449/8648640x^5 x^4$ &\\
$10, 12, 13]$&&$+ 1073561/1729728- 3994843/1081080x^3 + 17033563/2162160x^2 $&\\
\hline
\end{tabular}
\end{center}
}

{\small
\begin{center}
\begin{tabular}{|l|c|c|c|}
\hline
Orbit & result & Interpolated Map  & conjugation \\
\hline
\hline
\multicolumn{4}{|c|}{Many Preperiodic Points}\\
\hline
$[1, 2, -8, -1, -4, -6, -2, -7]$ & $13 (9, 1)$& $1/1260x^7 + 1/72x^6 + 23/360x^5 - 1/9x^4$&\\
&&$- 511/360x^3 - 137/72x^2 + 1829/420x + 1$&\\
\cline{2-4}
  &reduced& $1/1260x^7 - 1/72x^6 + 23/360x^5 + 1/36x^4$ & $\sm{1,-5,0,1}$ \\
&form&$- 211/360x^3 + 71/72x^2 - 207/140x$&\\
\hline
$[-1, -6, -2, -8, -9, -4, -5, -3]$ & $12 (6, 4)$ & $-1/1260x^7 - 1/45x^6 - 89/360x^5 - 25/18x^4$&\\
&&$- 1409/360x^3 - 323/90x^2 + 437/105x - 1$&\\
\hline
\end{tabular}
\end{center}
}

{\small
\begin{center}
\begin{tabular}{|l|c|c|c|}
\hline
Orbit & result & Interpolated Map &  conjugation \\

\hline
\hline
\multicolumn{4}{|c|}{Long Periodic Cycles}\\
\hline
\hline
$[-1, 1, -2, -6, -8, -5, -4, -3]$ & $(0, 10)$& $17/2520x^7 + 29/180x^6 + 509/360x^5 + 193/36x^4 $&\\
&&$+ 611/90x^3 - 226/45x^2 - 2039/210x - 1$ &\\
\cline{2-4}
 &reduced form& $17/2520x^7 - 1/36x^6 - 67/360x^5 + 23/36x^4$ & $\sm{1,-4,0,1}$ \\
&&$+ 131/90x^3 - 28/9x^2 - 583/210x + 1$&\\
\hline
%$[-1, -2, 2, -5, -4, -6, -3, 1]$ &$(0, 10)$ & $17/2520x^7 + 3/40x^6 + 11/90x^5 - 11/12x^4 - 781/360x^3 + 281/120x^2 + 1907/420x - 1$ &&\\
%reduced form && $17/2520x^7 + 1/36x^6 - 67/360x^5 - 23/36x^4 + 131/90x^3 + 28/9x^2 - 583/210x - 1$ & $[1,-1,0,1]$ &\\
%\hline
$[1, -7, -5, -4, -8, -1, -2, -3]$ & $(1, 9)$& $-11/2520x^7 - 73/720x^6 - 125/144x^5 - 479/144x^4$ &\\
&&$- 3889/720x^3 - 373/180x^2 + 317/84x + 1$&\\
\cline{2-4}
 &reduced form& unknown &\\
\hline
$[1, -1, -2, -6, -3, -4, 2, -5]$ & $(0, 9)$& $-19/5040x^7 - 3/80x^6 - 19/720x^5 + 23/48x^4$&\\
&&$+ 23/180x^3 - 353/120x^2 + 169/420x + 1$&\\
\cline{2-4}
 &reduced form& unknown&\\
\hline
\end{tabular}
\end{center}
}

{\small
\begin{center}
\begin{tabular}{|l|c|c|c|}
\hline
Orbit & result & Interpolated Map &  conjugation \\
\hline
\hline
\multicolumn{4}{|c|}{Long Preperiodic Tail}\\
\hline
$[1, 3, 5, 4, 6, -1, -3, -4]$ & $(9, 1)$& $11/20160x^7 - 41/10080x^6 - 173/10080x^5 + 10/63x^4$&\\
&&$- 1321/20160x^3 - 11639/10080x^2 + 863/280x + 1$&\\
\cline{2-4}
 &reduced & $11/315x^7 - 1/126x^6 - 607/1260x^5 + 125/504x^4$ & $\sm{2,1,0,1}$ \\
&form&$+ 1691/1260x^3 - 625/504x^2 + 116/105x + 1$&\\
\hline
$[1, 2, -8, -1, -4, -6, -2, -7]$ & $(9, 1)$& $1/1260x^7 + 1/72x^6 + 23/360x^5 - 1/9x^4$&\\
&&$- 511/360x^3 - 137/72x^2 + 1829/420x + 1$&\\
\cline{2-4}
 &reduced & $1/1260x^7 - 1/72x^6 + 23/360x^5 + 1/36x^4$ & $\sm{1,-5,0,1}$ \\
&form&$- 211/360x^3 + 71/72x^2 - 207/140x$&\\
\hline
$[-1, -7, -2, -5, -9,$ & $(8, 2)$& $-13/20160x^7 - 397/20160x^6 - 4171/20160x^5 - 3131/4032x^4$&\\
$ -4, -6, -8]$ &&$+ 961/2520x^3 + 40553/5040x^2 + 1373/105x - 1$&\\
\cline{2-4}
 &reduced & $-13/315x^7 - 11/210x^6 + 989/1260x^5 + 107/168x^4 $ & $\sm{2,-4,0,1}$ \\
&form&$- 2189/630x^3 - 701/840x^2 + 731/210x - 1$&\\
\hline
$[1, -6, 2, -2, -8,$ & $(8, 2)$& $1/252x^7 + 29/360x^6 + 49/90x^5 + 79/72x^4$&\\
$ -3, -5, -1]$&&$- 71/36x^3 - 323/45x^2 + 89/210x + 1$&\\
\cline{2-4}
 &reduced & $1/252x^7 - 11/360x^6 - 1/18x^5 + 47/72x^4$ & $\sm{1,-4,0,1}$ \\
&form&$+ 1/36x^3 - 118/45x^2 - 83/42x + 3$&\\
\hline
\end{tabular}
\end{center}
}

\begin{comment}
\subsubsection{Best Previously Known}
\begin{center}
\begin{tabular}{|l|c|c|}
\hline
Map & value & citation\\
\hline
\multicolumn{3}{|c|}{Long Periodic Cycles}\\
\hline
&&\\
\hline
\multicolumn{3}{|c|}{Long Preperiodic Tail}\\
\hline
&&\\\hline
\multicolumn{3}{|c|}{Many Preperiodic Points}\\
\hline
$1/5040x^7 - 1/90x^6 + 89/360x^5 - 25/9x^4 + 11947/720x^3 - 4609/90x^2 + 31147/420x - 30$& 16 (2,4) &Doyle-Hyde\\
\hline
\multicolumn{3}{|c|}{Small Height ratio}\\
\hline
&&\\
\hline
\end{tabular}
\end{center}
\end{comment}

%%%%%%%%%%%%%%%%%%%%%%%%%%%%%%%%%%%%%
\subsection{Degree 8: Polynomial}
{\small
\begin{center}
\begin{tabular}{|l|c|c|c|}
\hline
Orbit & result & Interpolated Map & conjugation\\
\hline
\hline
\multicolumn{4}{|c|}{Small Height Ratio}\\
\hline
$[1, 3, -10, -12, 2, -1, $ & $8.7710 \cdot 10^{-13}$ & $-1/39312x^8 - 71/83160x^7 - 265/30888x^6 $ &\\
$-9, -11, -8]$&&$- 571/77220x^5 + 1105/4752x^4 + 1549/154440x^3 $&\\
&&$- 304427/72072x^2 + 360253/60060x + 1$&\\
\cline{2-4}
 &reduced form& $-1/39312x^8 + 151/154440x^7 - 241/19305x^6 $ & $\sm{1,-9,0,1}$ \\
 &&$+ 6499/154440x^5 + 14995/61776x^4 - 1789/1485x^3$&\\
 &&$ - 623503/360360x^2 + 5677/8580x - 2$&\\
\hline
\end{tabular}
\end{center}
}
{\small
\begin{center}
\begin{tabular}{|l|c|c|c|}
\hline
Orbit & result & Interpolated Map & conjugation\\
\hline
\hline
\multicolumn{4}{|c|}{Many Preperiodic Points}\\
\hline
$[-1, 3, 2, -2, 1, -3,$ & $12 $ & $1/1260x^8 + 53/5040x^7 + 7/360x^6 $ &\\
$-6, -4, -7]$ && $- 17/90x^5 - 73/180x^4 + 1061/720x^3 $&\\
&&$+ 3491/2520x^2 - 451/105x - 1$&\\
\cline{2-4}
 &reduced form& $1/1260x^8 + 1/240x^7 - 23/720x^6 $ & $\sm{1,-1,0,1}$ \\
 &&$- 31/240x^5 + 373/720x^4 + 137/120x^3$&\\
 &&$ - 2929/840x^2 - 121/60x + 4$&\\
\hline
\end{tabular}
\end{center}
}
{\small
\begin{center}
\begin{tabular}{|l|c|c|c|}
\hline
Orbit & result & Interpolated Map & conjugation \\
\hline
\hline
\multicolumn{4}{|c|}{Long Periodic Cycles}\\
\hline
\hline
$[1, 7, 5, 8, 6, 9, 4, 2, 3]$ & $(0, 11)$ & $-1/8064x^8 + 1/672x^7 + 13/576x^6 - 29/60x^5$&\\
&&$+ 3377/1152x^4 - 595/96x^3 - 305/672x^2$ &\\
&&$+ 8551/840x + 1$&\\
\cline{2-4}
 &reduced form& unknown &\\
\hline
$[1, 9, 4, 8, 6, 2, 7, 5, 3]$ &$(0, 11)$& $-1/1920x^8 + 5/224x^7 - 1127/2880x^6 + 18/5x^5$&\\
&&$- 106747/5760x^4 + 5059/96x^3 - 110269/1440x^2$&\\
\hline
\end{tabular}
\end{center}
}{\small
\begin{center}
\begin{tabular}{|l|c|c|c|}
\hline
Orbit & result & Interpolated Map & conjugation\\
\hline
\hline
\multicolumn{4}{|c|}{Long Preperiodic Tail}\\
\hline
$[-1, -5, -7, -3, -8, -2, $; & $(9, 2)$ & $-1/4480x^8 - 73/10080x^7 - 259/2880x^6 - 187/360x^5 $&\\
$-6, -4, -9]$&&$- 7223/5760x^4 - 1/1440x^3 + 48829/10080x^2 $&\\
&&$+ 6743/840x - 1$&\\
\cline{2-4}
 &reduced form& $-1/4480x^8 + 1/288x^7 - 31/2880x^6 $ & $\sm{1,-6,0,1}$ \\
 &&$- 41/720x^5 + 1537/5760x^4 + 85/288x^3$&\\
 &&$ - 12659/10080x^2 - 269/120x + 2$&\\
\hline
\hline
\end{tabular}
\end{center}
}
\begin{comment}
\subsubsection{Best Previously Known}
\begin{center}
\begin{tabular}{|l|c|c|}
\hline
Map & value & citation\\
\hline
\multicolumn{3}{|c|}{Long Periodic Cycles}\\
\hline
&&\\
\hline
\multicolumn{3}{|c|}{Long Preperiodic Tail}\\
\hline
&&\\\hline
\multicolumn{3}{|c|}{Many Preperiodic Points}\\
\hline
$1/20160x^8 - 17/5040x^7 + 137/1440x^6 - 1037/720x^5 + 36427/2880x^4 - 2924/45x^3 + 313751/1680x^2 - 37111/140x + 142$& 17 (2,3) &Doyle-Hyde\\
\hline
\multicolumn{3}{|c|}{Small Height ratio}\\
\hline
&&\\
\hline
\end{tabular}
\end{center}
\end{comment}

%%%%%%%%%%%%%%%%%%%%%%%%%%%%%%%%%%%%%
\subsection{Degree 9: Polynomial}
{\small
\begin{center}
\begin{tabular}{|l|c|c|c|}
\hline
Orbit & result & Interpolated Map & conjugation \\
\hline
\hline
\multicolumn{4}{|c|}{Small Height Ratio}\\
\hline
$[-1, -12, -9, 1,   $ & $1.7173\cdot 10^{-14}$ & $2393/129729600x^9 + 13421/14414400x^8 + 49339/2702700x^7 $&\\
$-4,-3,-8, $&&$+ 563/3300x^6 + 334679/475200x^5 + 17903/52800x^4 $&\\
$-11, -10, 2]$&&$- 20931181/4054050x^3 - 27065231/3603600x^2 + 1520803/180180x - 1$&\\
\cline{2-4}
 &reduced form& $2393/129729600x^9 + 1283/4804800x^8 - 353/386100x^7$ & $\sm{1,-4,0,1}$ \\
 &&$- 323/14300x^6 + 9707/6177600x^5 + 139917/228800x^4$&\\
 &&$ + 181343/737100x^3 - 6346253/1201200x^2 - 2845/5148x + 1$&\\
\hline
 \hline
\end{tabular}
\end{center}
}

{\small
\begin{center}
\begin{tabular}{|l|c|c|c|}
\hline
Orbit & result & Interpolated Map & conjugation\\
\hline
\hline
\multicolumn{4}{|c|}{Long Preperiodic Tail}\\
\hline
$[1, -1, 6, 11, 8, 9, $ & $(10, 1)$& $19/332640x^9 - 839/332640x^8 + 15007/332640x^7$&\\
$10, 4, 3, 2]$&&$ - 138371/332640x^6 + 85733/41580x^5 - 835573/166320x^4$&\\
&&$ + 187673/55440x^3 + 192443/27720x^2 - 3462/385x + 1$&\\
\cline{2-4}
 &reduced form& $19/332640x^9 - 163/166320x^8 + 1027/332640x^7 $ & $\sm{1,3,0,1}$\\
 &&$+ 211/8316x^6 - 8629/66528x^5 - 14261/83160x^4 + 17663/15840x^3$&\\
 &&$ + 1025/1584x^2 - 2257/4620x - 1$&\\
\hline
\end{tabular}
\end{center}
}

{\small
\begin{center}
\begin{tabular}{|l|c|c|c|}
\hline
Orbit & result & Interpolated Map &  conjugation \\

\hline
\hline
\multicolumn{4}{|c|}{Long Periodic Cycles}\\
\hline
\hline
$[1, 3, 2, 4, 8, 9, $ & $(1, 11)$& $-17/120960x^9 + 29/4032x^8 - 1033/6720x^7$ &\\
$10, 6, 7, 5]$&&$+ 2557/1440x^6 - 69013/5760x^5 + 27415/576x^4 $&\\
&&$- 3218629/30240x^3 + 600373/5040x^2 - 671/14x + 1$&\\
\cline{2-4}
 &reduced form& $-34/945x^9 - 16/315x^8 + 3/5x^7 + 28/45x^6 $ & $\sm{2,6,0,1}$\\
 &&$- 1261/360x^5 - 107/45x^4 + 60269/7560x^3 $&\\
 &&$+ 3221/1260x^2 - 287/60x + 1/2$&\\
\hline
\end{tabular}
\end{center}
}

\begin{comment}

\subsubsection{Best Previously Known}
\begin{center}
\begin{tabular}{|l|c|c|}
\hline
Map & value & citation\\
\hline
\multicolumn{3}{|c|}{Long Periodic Cycles}\\
\hline
&&\\
\hline
\multicolumn{3}{|c|}{Long Preperiodic Tail}\\
\hline
&&\\\hline
\multicolumn{3}{|c|}{Many Preperiodic Points}\\
\hline
$1/181440x^9 - 1/2016x^8 + 571/30240x^7 - 19/48x^6 + 43069/8640x^5 - 1235/32x^4 + 8188031/45360x^3 - 242177/504x^2 + 40382/63x - 306$& 20 (4,1) &Doyle-Hyde\\
\hline
\multicolumn{3}{|c|}{Small Height ratio}\\
\hline
&&\\
\hline
\end{tabular}
\end{center}
\end{comment}

%%%%%%%%%%%%%%%%%%%%%%%%%%%%%%%%%%%%%
\subsection{Degree 10: Polynomial}
{\small
\begin{center}
\begin{tabular}{|l|c|c|c|}
\hline
Orbit & result & Interpolated Map  & conjugation \\
\hline\hline
\multicolumn{4}{|c|}{Small Height ratio}\\
\hline
%$[1, 9, 8, 13, 7, 12, -2, 10, 2, 11, -1]$ & $3.1858 \cdot 10^{-15}$& $-23/103783680*x^{10} + 43/1886976*x^9 - 14213/17297280*x^8 + 247927/17297280*x^7 - 11645/89856*x^6 + 2706217/4942080*x^5 - 9766333/25945920*x^4 - 98357593/25945920*x^3 + 2164621/360360*x^2 + 516119/90090*x + 1$&&\\
%reduced form && $-23/103783680*x^{10} + 151/20756736*x^9 + 197/1572480*x^8 - 10471/17297280*x^7 - 48443/4942080*x^6 + 127319/4942080*x^5 + 7463537/25945920*x^4 - 15910571/25945920*x^3 - 388711/120120*x^2 + 204779/45045*x + 5$ & $[1,7,0,1]$ & $z=-7$\\
%\hline
%$[1, 12, 14, 9, 2, 10, 8, 3, -1, -2, 13]$ & $5.5071\cdot 10^{-15}$& $-1/982800x^{10} + 161/2316600x^9 - 20597/10810800x^8 + 143449/5405400x^7 - 137873/720720x^6 + 386902/675675x^5 + 5818997/10810800x^4 - 24363293/4054050x^3 + 19756453/5405400x^2 + 447221/36036x + 1$ &&\\
%reduced form && $-1/982800x^{10} - 1/579150x^9 + 1241/5405400x^8 + 233/415800x^7 - 168071/10810800x^6 - 10273/300300x^5 + 63983/163800x^4 + 677521/1247400x^3 - 5327293/1801800x^2 + 21767/180180x - 293/143$ & $z=-7$\\
%\hline
%$[1, 8, 10, -2, 9, 3, 14, -1, 13, 12, 11];9.29e-15;0$&&&&&\\
%$[1, 9, 11, 2, 4, 12, -1, 13, 14, 10, 3];1.55155376266151e-15;0$&&&&&\\
$[-1, 3, 15, 9, 12, 1, $ & $9.9368\cdot 10^{-16}$ &\scalebox{0.8}{$947/518918400x^{10} - 64921/518918400x^9 + 11029/3088800x^8 $}&\\
$10, 11, 4, 14, 2]$&&\scalebox{0.8}{$- 431149/7862400x^7 + 28151/57600x^6 - 1608131/633600x^5 $}&\\
&&\scalebox{0.8}{$+ 1885570259/259459200x^4 - 1239296549/129729600x^3 $}&\\
&&\scalebox{0.8}{$- 267299/1029600x^2 + 5638181/360360x - 1$}&\\
\cline{2-4}
 &reduce form& \scalebox{0.8}{$947/518918400x^{10} + 20309/518918400x^9 + 7681/86486400x^8 $} & $\sm{1,9,0,1}$ \\
 &&\scalebox{0.8}{$- 252521/86486400x^7 - 9167/633600x^6 + 47239/633600x^5 $}&\\
\hline
\end{tabular}
\end{center}
}
{\small
\begin{center}
\begin{tabular}{|l|c|c|c|}
\hline
Orbit & result & Interpolated Map & conjugation \\
\hline
\hline
\multicolumn{4}{|c|}{Many Preperiodic Points}\\
\hline
\hline
\end{tabular}
\end{center}
}

{\small
\begin{center}
\begin{tabular}{|l|c|c|c|}
\hline
Orbit & result & Interpolated Map & conjugation \\
\hline
\hline
\multicolumn{4}{|c|}{Long Periodic Cycles}\\
\hline
\hline
$[1, 2, 3, -6, -11, -2, $ & $(0, 12)$& \scalebox{0.8}{$-37/64864800x^{10} + 17/9266400x^9 + 17467/43243200x^8 $} &\\
$-7, 4, -9, -10, -1]$&&\scalebox{0.8}{$+ 781/163800x^7+ 347/64350x^6 - 411569/3088800x^5 $}&\\
&&\scalebox{0.8}{$ - 37024291/129729600x^4 + 4938731/4633200x^3$}&\\
&&\scalebox{0.8}{$ + 91597/327600x^2+ 11267/180180x + 1$}&\\
\cline{2-4}
 &reduced form& same &$\sm{1,0,0,1}$\\
\hline
%$[1, -12, -10, -13, -4, -11, -2, -3, -8, -6, -1];-60;(0, 12)$&&&&&\\
%$[1, -10, -4, -1, -12, 3, -11, 2, -8, -9, -7];-60;(0, 12)$&&&&&\\
%$[-1, -11, -6, -2, -4, -10, 1, -9, -7, -8, -5];-60;(0, 12)$&&&&&\\
%$[-1, -7, -9, -6, -3, -4, 2, -12, -10, 1, -11];-60;(0, 12)$&&&&&\\
$[-1, 1, 2, 3, -12, -3, $ & $(0, 12)$& \scalebox{0.8}{$3607/1556755200x^{10} + 1579/14826240x^9 + 23641/12972960x^8 $}&\\
\hline
\end{tabular}
\end{center}
}
{\small
\begin{center}
\begin{tabular}{|l|c|c|c|}
\hline
Orbit & result & Interpolated Map  & conjugation \\
\hline
\hline
\multicolumn{4}{|c|}{Long Preperiodic Tail}\\
\hline
$[-1, 2, 4, -7, -5, 1,$ & $(11, 1)$& \scalebox{0.8}{$53/6652800x^{10} + 41/133056x^9 + 2351/665280x^8 $}&\\
$ -4, -3, -8, 3, -6]$&&\scalebox{0.8}{$+ 2731/332640x^7 - 574771/6652800x^6 - 269707/665280x^5$}&\\ 
&&\scalebox{0.8}{$+ 25087/55440x^4 + 288641/83160x^3$}&\\
&&\scalebox{0.8}{$- 8539/23100x^2 - 4677/770x - 1$}&\\
\cline{2-4}
 &reduced form& \scalebox{0.8}{$212/51975x^{10} + 4/105x^9 - 2/27x^8 $} & $\sm{2,-2,0,1}$ \\
 &&\scalebox{0.8}{$- 26/35x^7 + 317/700x^6 + 283/60x^5 $}&\\
\hline
\end{tabular}
\end{center}
}

%%%%%%%%%%%%%%%%%%%%%%%%%%%%%%%%%%%%%
\subsection{Degree 11: Polynomial}
{\small
\begin{center}
\begin{tabular}{|l||c|c|c|}
\hline
Orbit & result & Interpolated Map  & conjugation\\
\hline
\hline\multicolumn{4}{|c|}{Small Height ratio}\\
\hline
$[-1, 1, -10, $ & $7.3446\cdot 10^{-17}$ & \scalebox{0.8}{$47/172972800x^{11} + 439/24710400x^{10} + 12553/25945920x^9 $} &\\
$-9, -12, -11,$&&\scalebox{0.8}{$+ 60947/8648640x^8 + 10051411/172972800x^7 + 234841/915200x^6 $}&\\
$-8, -14, -3,$&&\scalebox{0.8}{$+ 190819/443520x^5 - 143957/224640x^4 - 91763489/32432400x^3 $}&\\
$ -2, 2, -13]$&&\scalebox{0.8}{$- 5193893/1663200x^2 - 1138541/360360x - 1$}&\\
\cline{2-4}
 &reduced form& \scalebox{0.8}{$47/172972800x^{11} + 71/4804800x^{10} + 6665/20756736x^9 $} & $\sm{1,-1,0,1}$ \\
 &&\scalebox{0.8}{$+ 59629/17297280x^8 + 2959361/172972800x^7 + 15619/1544400x^6 $}&\\
\hline
\end{tabular}
\end{center}
}
{\small
\begin{center}
\begin{tabular}{|l|c|c|c|}
\hline
\hline
\multicolumn{4}{|c|}{Many PrePeriodic Points}\\
\hline
Orbit & result & Interpolated Map & conjugation\\
\hline
\end{tabular}
\end{center}
}

{\small
\begin{center}
\begin{tabular}{|l|c|c|c|}
\hline
\hline
\multicolumn{4}{|c|}{Long Periodic Cycles}\\
\hline
Orbit & result & Interpolated Map & conjugation\\
\hline
\hline
%$[1, 3, 14, 6, 9, -2, 8, 2, 11, -1, 13, 12];-65;(0, 13)$&&&&&\\
$[1, 7, 4, -2, 11, 6, -4,$ & $(0, 13)$& \scalebox{0.8}{$-257/518918400x^{11} + 2657/103783680x^{10} - 27829/51891840x^9 $}&\\
$ -1, 12, 9, 8, 10]$&&\scalebox{0.8}{$+ 48469/8648640x^8 - 4519987/172972800x^7 - 80383/4942080x^6 $}&\\
&&\scalebox{0.8}{$+ 17404669/25945920x^5 - 100793989/51891840x^4 - 112467679/64864800x^3 $}&\\
&&\scalebox{0.8}{$+ 1159031/110880x^2 - 508201/360360x + 1$}&\\
\cline{2-4}
 &reduced form& \scalebox{0.8}{$-257/518918400x^{11} - 1333/74131200x^{10} - 12013/51891840x^9 $} & $\sm{1,8,0,1}$ \\
 &&\scalebox{0.8}{$- 1931/1729728x^8 + 76493/172972800x^7 + 544877/24710400x^6 $}&\\
\end{tabular}
\end{center}
}
{\small
\begin{center}
\begin{tabular}{|l|c|c|c|}
\hline
Orbit & result & Interpolated Map  & conjugation \\
\hline
\hline
\multicolumn{4}{|c|}{Long Preperiodic Tail}\\
\hline
$[-1, -5, -7, $ & $(12, 1)$& \scalebox{0.8}{$-173/43243200x^{11} - 29671/172972800x^{10} - 59779/20756736x^9 $}&\\
$-2, -8, -11,$&&\scalebox{0.8}{$- 397403/17297280x^8 - 6400363/86486400x^7 + 1753/19200x^6 $}&\\
$3,-6, 1,$&&\scalebox{0.8}{$+ 41746499/34594560x^5 + 17647919/8648640x^4 - 332984849/129729600x^3 $}&\\
$ -10, -9, 2]$&&\scalebox{0.8}{$- 23266729/2702700x^2 - 21283/20020x - 1$}&\\
\cline{2-4}
 &reduced form& \scalebox{0.8}{$-173/43243200x^{11} - 817/6406400x^{10} - 28741/20756736x^9 $} & $\sm{1,-1,0,1}$ \\
 &&\scalebox{0.8}{$- 391/95040x^8 + 2194987/86486400x^7 + 4294159/24710400x^6 $}&\\
\hline
\end{tabular}
\end{center}
}

%%%%%%%%%%%%%%%%%%%%%%%%%%%%%%%%%%%%%
\subsection{Degree 12: Polynomial}
{\small
\begin{center}
\begin{tabular}{|l|c||c|c|}
\hline
Orbit & result & Interpolated Map  & conjugation\\
\hline
\hline
\multicolumn{4}{|c|}{Small Height ratio}\\
\hline
$[1, -2, -11, -8, $ & $5.9715 \cdot 10^{-19}$ & \scalebox{0.8}{$1592279639/9714712879872000x^{12} + 6499459/19989121152000x^{11} $}&\\
$-1, 12, 7, $&&\scalebox{0.8}{$- 642202438277/9714712879872000x^{10} + 221170236173/1619118813312000x^9  $}&\\
$3, -18, 8,$&&\scalebox{0.8}{$+ 4948365818327/647647525324800x^8 - 166647918583/4693098009600x^7 $}&\\
$ 5, 2, -15]$&&\scalebox{0.8}{$- 332836924433609/1387816125696000x^6 + 44212194839429/25700298624000x^5 $}&\\ 
&&\scalebox{0.8}{$- 3535545209345011/2428678219968000x^4- 2375054642323511/404779703328000x^3 $}&\\
&&\scalebox{0.8}{$+ 38373129132491/6746328388800x^2- 2514169733/892371480x + 1$}&\\
\hline
\hline
\end{tabular}
\end{center}
}

{\small
\begin{center}
\begin{tabular}{|l|c|c|c|}
\hline
Orbit & result & Interpolated Map & conjugation \\
\hline
\hline
\multicolumn{4}{|c|}{Many Preperiodic Points}\\
\hline
\end{tabular}
\end{center}
}

{\small
\begin{center}
\begin{tabular}{|l|c|c|c|}
\hline
Orbit & result & Interpolated Map & conjugation \\
\hline
\hline
\multicolumn{4}{|c|}{Long Periodic Cycles}\\
\hline
\hline
$[-1, -6, 7, -7, $ & $(0, 14)$& \scalebox{0.8}{$-6847/21794572800x^{12} - 109/28828800x^{11} + 743473/21794572800x^{10}$}&\\
$-8, 4, 6, $&&\scalebox{0.8}{$ + 5003/10378368x^9 - 8656007/7264857600x^8 - 1906781/86486400x^7 $}&\\
$-3, 5, -4, $&&\scalebox{0.8}{$+ 283709539/21794572800x^6 + 7538287/17297280x^5 + 39484223/1362160800x^4$}&\\
$1, 3, -5]$&&\scalebox{0.8}{$ - 113189047/32432400x^3 - 81857723/151351200x^2 $}&\\
&&\scalebox{0.8}{$+ 2730001/360360x - 1$}&\\
\cline{2-4}
 &reduced form& \scalebox{0.8}{$-6847/21794572800x^{12} + 20521/1816214400x^{11} - 220091/1676505600x^{10} $}& $\sm{1,-4,0,1}$ \\
\hline
\end{tabular}
\end{center}
}
{\small
\begin{center}
\begin{tabular}{|l|c|c|c|}
\hline
Orbit & result & Interpolated Map  & conjugation \\
\hline
\hline
\multicolumn{4}{|c|}{Long Preperiodic Tail}\\
\hline
$[1, 10, 13, 3,  $ & $(13, 1)$ & \scalebox{0.8}{$-16391/43589145600x^{12} + 177811/7264857600x^{11} - 5911127/8717829120x^{10} $} &\\
$6, 9, 2,$&&\scalebox{0.8}{$+ 15085313/1452971520x^9 - 1379902411/14529715200x^8 + 416084377/807206400x^7 $}&\\
$5, 12, -1,$&&\scalebox{0.8}{$- 12679380821/8717829120x^6 + 945170263/1452971520x^5  $}&\\
$ 4, 11, -2]$&&\scalebox{0.8}{$+ 80268528481/10897286400x^4- 27942029851/1816214400x^3 $}&\\ 
&&\scalebox{0.8}{$+ 2787133/15135120x^2+ 43408499/2522520x + 1$}&\\
\cline{2-4}
 &reduced form& \scalebox{0.8}{$-16391/43589145600x^{12} + 112247/7264857600x^{11} - 1487401/6227020800x^{10} $} & $\sm{1,2,0,1}$ \\
\hline
\end{tabular}
\end{center}
}

{\small
\begin{center}
\begin{tabular}{|l|c||c|c|}
\hline
Orbit & result & Interpolated Map  & conjugation\\
\hline
\hline
\multicolumn{4}{|c|}{Small Height ratio}\\
\hline
$[1, -2, -11, -8, $ & $5.9715 \cdot 10^{-19}$ & \scalebox{0.8}{$1592279639/9714712879872000x^{12} + 6499459/19989121152000x^{11}  $}&\\
$-1, 12, 7,$&&\scalebox{0.8}{$- 642202438277/9714712879872000x^{10}+ 221170236173/1619118813312000x^9  $}&\\
$3, -18, 8,$&&\scalebox{0.8}{$+ 4948365818327/647647525324800x^8 - 166647918583/4693098009600x^7$}&\\
$ 5, 2, -15]$&&\scalebox{0.8}{$- 332836924433609/1387816125696000x^6 + 44212194839429/25700298624000x^5 $}&\\
&&\scalebox{0.8}{$ - 3535545209345011/2428678219968000x^4- 2375054642323511/404779703328000x^3 $}&\\ 
&&\scalebox{0.8}{$+ 38373129132491/6746328388800x^2- 2514169733/892371480x + 1$}&\\
\hline
\end{tabular}
\end{center}
}

%%%%%%%%%%%%%%%%%%%%%%%%%%%%%%%%%%%%%
\subsection{Degree 13: Polynomial}
{\small
\begin{center}
\begin{tabular}{|l|c|c|c|}
\hline
Orbit & result & Interpolated Map & conjugation \\
\hline
\hline
\multicolumn{4}{|c|}{Small Height Ratio}\\
\hline
\end{tabular}
\end{center}
}
{\small
\begin{center}
\begin{tabular}{|l|c|c|c|}
\hline
Orbit & result & Interpolated Map & conjugation \\
\hline
\hline
\multicolumn{4}{|c|}{Many Preperiodic Points}\\
\hline
\end{tabular}
\end{center}
}
{\small
\begin{center}
\begin{tabular}{|l|c|c|c|}
\hline
Orbit & result & Interpolated Map  & conjugation \\
\hline
\hline
\multicolumn{4}{|c|}{Long Periodic Cycles}\\
\hline
\hline
$[-1, -3, -2, -13,  $ & $(0, 15)$& \scalebox{0.7}{$-43/311351040x^{13} - 949/95800320x^{12} - 73313/239500800x^{11} $}&\\
$1, -12, -10,-11, $&&\scalebox{0.7}{$- 230831/43545600x^{10} - 16183/290304x^9 - 1032019/2903040x^8 $}&\\
$-5, -7, -8,$&&\scalebox{0.7}{$- 26689249/21772800x^7 - 40975853/43545600x^6 + 36782723/4354560x^5 $}&\\
$ -6, -9, -4]$&&\scalebox{0.7}{$+ 60055129/2177280x^4 + 115981307/6652800x^3 - 27260551/831600x^2$}&\\ 
&&\scalebox{0.7}{$-499331/17160x - 1$}&\\
\cline{2-4}
 &reduced form& unknown&\\
\hline
\end{tabular}
\end{center}
}
{
\small
\begin{center}
\begin{tabular}{|l|c|c|c|}
\hline
Orbit & result & Interpolated Map  & conjugation \\
\hline
\hline
\multicolumn{4}{|c|}{Long Preperiodic Tail}\\
\hline
$[-1, -2, -4, -12,  $ & $(14, 1)$& \scalebox{0.7}{$4177/21794572800x^{13} + 699763/43589145600x^{12} + 584513/990662400x^{11} $}&\\
$-14, -8, -10,-3, $&&\scalebox{0.7}{$+ 21782699/1743565824x^{10} + 12724267/75675600x^9 + 21727697183/14529715200x^8 $}&\\ 
$-5, -7, -13,$&&\scalebox{0.7}{$+ 96254583559/10897286400x^7 + 58839678665/1743565824x^6 + 1680904706167/21794572800x^5 $}&\\ 
$ -11, 1, -6]$&&\scalebox{0.7}{$+ 896054795047/10897286400x^4 - 35559014017/1816214400x^3 - 364701871/3027024x^2$}&\\ 
&&\scalebox{0.7}{$ - 24701923/360360x - 1$}&\\ 
\cline{2-4}
 &reduced form& \scalebox{0.7}{$4177/21794572800x^{13} + 48151/43589145600x^{12} - 301583/10897286400x^{11} $} & $\sm{1,-6,0,1}$\\
\hline

\end{tabular}
\end{center}
}

\end{all_data}

\end{document}